# ON THE ERGODICITY PROPERTIES OF SOME ADAPTIVE MCMC ALGORITHMS

By Christophe Andrieu and Éric Moulines

*University of Bristol and École Nationale Supérieure des Télécommunications*

In this paper we study the ergodicity properties of some adaptive Markov chain Monte Carlo algorithms (MCMC) that have been recently proposed in the literature. We prove that under a set of verifiable conditions, ergodic averages calculated from the output of a so-called adaptive MCMC sampler converge to the required value and can even, under more stringent assumptions, satisfy a central limit theorem. We prove that the conditions required are satisfied for the independent Metropolis–Hastings algorithm and the random walk Metropolis algorithm with symmetric increments. Finally, we propose an application of these results to the case where the proposal distribution of the Metropolis–Hastings update is a mixture of distributions from a curved exponential family.

**1. Introduction.** Markov chain Monte Carlo (MCMC), introduced by Metropolis et al. [27], is a popular computational method for generating samples from virtually any distribution $\pi$. In particular, there is no need for the normalizing constant to be known, and the space $\mathsf{X} \subset \mathbb{R}^{n_x}$ (for some integer $n_x$) on which it is defined can be high dimensional. We will hereafter denote by $\mathcal{B}(\mathsf{X})$ the associated countably-generated $\sigma$-field. The method consists of simulating an ergodic Markov chain $\{X_k, k \geq 0\}$ on $\mathsf{X}$ with transition probability $P$ such that $\pi$ is a *stationary* distribution for this chain, that is, $\pi P = \pi$. Such samples can be used, for example, to compute integrals of the form

$$\pi(f) \stackrel{\text{def}}{=} \int_\mathsf{X} f(x)\pi(dx),$$

for some $\pi$-integrable function $f: \mathsf{X} \to \mathbb{R}$, using estimators of the type

$$(1) \qquad S_n(f) = \frac{1}{n}\sum_{k=1}^{n} f(X_k).$$

In general, the transition probability $P$ of the Markov chain depends on some tuning parameter, say $\theta$, defined on some space $\Theta \subset \mathbb{R}^{n_\theta}$ for some integer $n_\theta$, and the convergence properties of the Monte Carlo averages in equation (1) might strongly depend on a proper choice of these parameters.

We illustrate this here with the Metropolis–Hastings (MH) update, but it should be stressed at this point that the results presented in this paper apply to much more general settings (including, in particular, hybrid samplers, sequential or population Monte Carlo samplers). The MH algorithm requires the choice of a *proposal distribution* $q$. In order to simplify the discussion, we will here assume that $\pi$ and $q$ admit densities with respect to the Lebesgue measure $\lambda^{\mathrm{Leb}}$, denoted (with an abuse of notation) $\pi$ and $q$ hereafter. The role of the distribution $q$ consists of proposing potential transitions for the Markov chain $\{X_k\}$. Given that the chain is currently at $x$, a candidate $y$ is accepted with probability $\alpha(x,y)$, defined as

$$\alpha(x,y) = \begin{cases} 1 \wedge \dfrac{\pi(y)}{\pi(x)}\dfrac{q(y,x)}{q(x,y)}, & \text{if } \pi(x)q(x,y) > 0, \\ 1, & \text{otherwise,} \end{cases}$$

where $a \wedge b \stackrel{\mathrm{def}}{=} \min(a,b)$. Otherwise it is rejected, and the Markov chain stays at its current location $x$. For $(x, A) \in \mathsf{X} \times \mathcal{B}(\mathsf{X})$, the transition kernel $P$ of this Markov chain takes the form

$$\begin{aligned}(2) \quad P(x, A) &= \int_{A-x} \alpha(x, x+z) q(x, x+z) \lambda^{\mathrm{Leb}}(dz) \\ &\quad + \mathbb{1}_A(x) \int_{\mathsf{X}-x} (1-\alpha(x, x+z)) q(x, x+z) \lambda^{\mathrm{Leb}}(dz),\end{aligned}$$

where $A - x \stackrel{\mathrm{def}}{=} \{z \in \mathsf{X}, x+z \in A\}$. The Markov chain $P$ is reversible with respect to $\pi$ and therefore admits $\pi$ as an invariant distribution. Conditions on the proposal distribution $q$ that guarantee irreducibility and positive recurrence are mild, and many satisfactory choices are possible. For the purpose of illustration, we concentrate in this introduction on the symmetric increments random walk MH algorithm (hereafter SRWM), in which $q(x,y) = q(y-x)$ for some symmetric probability density $q$ on $\mathbb{R}^{n_x}$, referred to as the *increment distribution*. The transition kernel $P_q^{\mathrm{SRW}}$ of the Metropolis algorithm is then given for $x, A \in \mathsf{X} \times \mathcal{B}(\mathsf{X})$ by

$$P_q^{\mathrm{SRW}}(x, A) = \int_{A-x} \left(1 \wedge \frac{\pi(x+z)}{\pi(x)}\right) q(z) \lambda^{\mathrm{Leb}}(dz)$$

(3)
$$+ \mathbb{1}_A(x) \int_{\mathsf{X}-x} \left(1 - 1 \wedge \frac{\pi(x+z)}{\pi(x)}\right) q(z) \lambda^{\text{Leb}}(dz).$$

A classical choice for $q$ is the multivariate normal distribution with zero mean and covariance matrix $\Gamma$, $\mathcal{N}(0,\Gamma)$. We will later on refer to this algorithm as the N-SRWM. It is well known that either too small or too large a covariance matrix will result in highly positively correlated Markov chains and therefore estimators $S_n(f)$ with a large variance. Gelman, Roberts and Gilks [16] have shown that the "optimal" covariance matrix (under restrictive technical conditions not given here) for the N-SRWM is $(2.38^2/n_x)\Gamma_\pi$, where $\Gamma_\pi$ is the true covariance matrix of the target distribution. In practice, this covariance matrix $\Gamma$ is determined by trial and error, using several realizations of the Markov chain. This hand-tuning requires some expertise and can be time consuming. In order to circumvent this problem, Haario, Saksman and Tamminen [19] have proposed to "learn $\Gamma$ on the fly." The Haario, Saksman and Tamminen [19] algorithm can be summarized as follows:

(4)
$$\mu_{k+1} = \mu_k + \gamma_{k+1}(X_{k+1} - \mu_k), \qquad k \geq 0,$$
$$\Gamma_{k+1} = \Gamma_k + \gamma_{k+1}((X_{k+1} - \mu_k)(X_{k+1} - \mu_k)^{\mathrm{T}} - \Gamma_k),$$

where, denoting by $\mathcal{C}_+^{n_x}$ the cone of positive $n_x \times n_x$ matrices:

- $X_{k+1}$ is drawn from $P_{\theta_k}(X_k,\cdot)$, where for $\theta = (\mu,\Gamma) \in \Theta = \mathbb{R}^{n_x} \times \mathcal{C}_+^{n_x}$, $P_\theta = P_{\mathcal{N}(0,\lambda\Gamma)}^{\text{SRW}}$ is here the kernel of a symmetric random walk MH with a Gaussian increment distribution $\mathcal{N}(0,\lambda\Gamma)$, $\lambda > 0$ being a constant scaling factor depending only on the dimension of the state space $n_x$ and kept constant across the iterations;
- $\{\gamma_k\}$ is a nonincreasing sequence of positive stepsizes such that $\sum_{k=1}^\infty \gamma_k = \infty$ and $\sum_{k=1}^\infty \gamma_k^{1+\delta} < \infty$ for some $\delta > 0$ (Haario, Saksman and Tamminen [19] have suggested the choice $\gamma_k = 1/k$).

It was realized in [4] that such a scheme is a particular case of a more general framework. More precisely, for $\theta = (\mu,\Gamma) \in \Theta$, define $H:\Theta \times \mathsf{X} \to \mathbb{R}^{n_x} \times \mathbb{R}^{n_x \times n_x}$ as

(5) $$H(\theta,x) \stackrel{\text{def}}{=} (x - \mu, (x-\mu)(x-\mu)^{\mathrm{T}} - \Gamma)^{\mathrm{T}}.$$

With this notation, the recursion in (4) may be written as

(6) $$\theta_{k+1} = \theta_k + \gamma_{k+1} H(\theta_k, X_{k+1}), \qquad k \geq 0,$$

with $X_{k+1} \sim P_{\theta_k}(X_k,\cdot)$. This recursion is at the core of most of classical stochastic approximation algorithms (see, e.g., [9, 13, 23] and the references therein). This algorithm is designed to solve the equations $h(\theta) = 0$ where $\theta \mapsto h(\theta)$ is the so-called *mean field* defined as

(7) $$h(\theta) \stackrel{\text{def}}{=} \int_\mathsf{X} H(\theta,x) \pi(dx).$$

For the present example, assuming that $\int_X |x|^2 \pi(dx) < \infty$, one can easily check that

$$(8) \quad h(\theta) = \int_X H(\theta, x)\pi(dx) = (\mu_\pi - \mu, (\mu_\pi - \mu)(\mu_\pi - \mu)^T + \Gamma_\pi - \Gamma)^T,$$

with $\mu_\pi$ and $\Gamma_\pi$ the mean and covariance of the target distribution, that is,

$$(9) \quad \mu_\pi = \int_X x\pi(dx) \quad \text{and} \quad \Gamma_\pi = \int_X (x - \mu_\pi)(x - \mu_\pi)^T \pi(dx).$$

One can rewrite (6) as

$$\theta_{k+1} = \theta_k + \gamma_{k+1} h(\theta_k) + \gamma_{k+1} \xi_{k+1},$$

where $\{\xi_k = H(\theta_{k-1}, X_k) - h(\theta_{k-1}), k \geq 1\}$ is generally referred to as "the noise." The general theory of *stochastic approximation* (SA) provides us with conditions under which this recursion eventually converges to the set $\{\theta \in \Theta, h(\theta) = 0\}$. These issues are discussed in Sections 3 and 5.

In the context of adaptive MCMC, the parameter convergence is not the central issue; the focus is rather on the approximation of $\pi(f)$ by the sample mean $S_n(f)$. However, there is here a difficulty with the adaptive approach: as the parameter estimate $\theta_k = \theta_k(X_0, \ldots, X_k)$ depends on the whole past, the successive draws $\{X_k\}$ do not define an homogeneous Markov chain, and standard arguments for the consistency and asymptotic normality of $S_n(f)$ do not apply in this framework. Note that this is despite the fact that, for any $\theta \in \Theta$, $\pi P_\theta = \pi$. This is illustrated by the following example. Let $X = \{1, 2\}$ and consider for $\theta, \theta(1), \theta(2) \in \Theta$ the following Markov transition probability matrices:

$$P_\theta = \begin{bmatrix} 1-\theta & \theta \\ \theta & 1-\theta \end{bmatrix}, \qquad \widetilde{P} = \begin{bmatrix} 1-\theta(1) & \theta(1) \\ \theta(2) & 1-\theta(2) \end{bmatrix}.$$

One can check that for any $\theta \in \Theta$, $\pi = (1/2, 1/2)$ satisfies $\pi P_\theta = \pi$. However, if we let $\theta_k$ be a given function $\theta : X \to (0, 1)$ of the current state, that is, $\theta_k = \theta(X_k)$, one defines a new Markov chain with transition probability $\widetilde{P}$ now having $[\theta(2)/(\theta(1) + \theta(2)), \theta(1)/(\theta(1) + \theta(2))]$ as invariant distribution. One recovers $\pi$ when the dependence on the current state $X_k$ is removed or vanishes with the iterations. With this example in mind, the problems that we address in the present paper and our main general results can be summarized as follows:

1. In situations where $|\theta_{k+1} - \theta_k| \to 0$ as $k \to +\infty$ w.p. 1, we prove a strong law of large numbers for $S_n(f)$ (see Theorem 8) under mild additional conditions. Such a consistency result may arise even in situations where the parameter sequence $\{\theta_k\}$ does not converge.
2. In situations where $\theta_k$ converges w.p. 1, we prove an invariance principle for $\sqrt{n}(S_n(f) - \pi(f))$; the limiting distribution is, in general, a mixture of Gaussian distributions (see Theorem 9).

Note that Haario, Saksman and Tamminen [19] have proved the consistency of Monte Carlo averages for the algorithm described by (4). Our results apply to more general settings and rely on assumptions which are less restrictive than those used in [19]. The second point above, the invariance principle, has, to the best of our knowledge, not been addressed for adaptive MCMC algorithms. We point out that Atchadé and Rosenthal [6] have independently extended the consistency result of Haario, Saksman and Tamminen [19] to the case where X is unbounded, using Haario, Saksman and Tamminen's [19] mixingale technique. Our technique of proof is different, and our algorithm allows for an unbounded parameter $\theta$ to be considered, as opposed to Atchadé and Rosenthal [6].

The paper is organized as follows. In Section 2 we detail our general procedure and introduce some notation. In Section 3 we establish the consistency (i.e., a strong law of large numbers) for $S_n(f)$ (Theorem 8). In Section 4 we strengthen the conditions required to ensure the law of large numbers (LLN) for $S_n(f)$ and establish an invariance principle (Theorem 9). In Section 5 we focus on the classical Robbins–Monro implementation of our procedure and introduce further conditions that allow us to prove that $\{\theta_k\}$ converges w.p. 1 (Theorem 11). In Section 6 we establish general properties of the generic SRWM required to ensure an LLN and an invariance principle. For pedagogical purposes we show how to apply these results to the N-SRWM of [19] (Theorem 15). In Section 7 we present another application of our theory. We focus on the independent Metropolis–Hastings algorithm (IMH) and establish general properties required for the LLN and the invariance principle. We then go on to propose and analyse an algorithm that matches the so-called proposal distribution of the IMH to the target distribution $\pi$, in the case where the proposal distribution is a mixture of distributions from the exponential family. The main result of this section is Theorem 21. We conclude with the remark that this latter result equally applies to a generalization of the N-SRWM where the proposal is again a mixture of distributions. Application to samplers which consist of a mixture of SRWM and IMH is straightforward.

**2. Algorithm description and main definitions.** Before describing the procedure under study, it is necessary to introduce some notation and definitions. Let T be a separable space and let $\mathcal{B}(\mathsf{T})$ be a countably-generated $\sigma$-field on T. For a Markov chain with transition probability $P:\mathsf{T}\times\mathcal{B}(\mathsf{T})\to[0,1]$ and any nonnegative measurable function $f:\mathsf{T}\to[0,+\infty)$, we define $Pf(t)=P(t,f)\stackrel{\mathrm{def}}{=}\int_\mathsf{T} P(t,dt')f(t')$ and for any integer $k$, denote by $P^k$ the $k$th iterate of the kernel. For a probability measure $\mu$, we define, for any $A\in\mathcal{B}(\mathsf{T})$, $\mu P(A)\stackrel{\mathrm{def}}{=}\int_\mathsf{T}\mu(dt)P(t,A)$. A Markov chain on a state space T is said to be $\mu$-*irreducible* if there exists a measure $\mu$ on $\mathcal{B}(\mathsf{T})$ such that, whenever $\mu(A) > 0$, $\sum_{k=0}^{\infty} P^k(t, A) > 0$ for all $t \in \mathsf{T}$. Denote by $\mu$ a maximal irreducibility measure for $P$ (see [28], Chapter 4, for the definition and the construction of such a measure). If $P$ is $\mu$-irreducible, aperiodic and has an invariant probability measure $\pi$, then $\pi$ is unique and is a maximal irreducibility measure.

Two main ingredients are required for the definition of our adaptive MCMC algorithms:

1. A family of Markov transition kernels on $\mathsf{X}$, $\{P_\theta, \theta \in \Theta\}$ indexed by a finite-dimensional parameter $\theta \in \Theta \subset \mathbb{R}^{n_\theta}$, where $\Theta$ is asumed to be an open set. For each $\theta$ in $\Theta$, it is assumed that $P_\theta$ is $\pi$-irreducible and that $\pi P_\theta = \pi$, that is, $\pi$ is the invariant distribution for $P_\theta$.
2. A family of *update functions* $\{H(\theta, x) : \Theta \times \mathsf{X} \mapsto \mathbb{R}^{n_\theta}\}$ which are used to adapt the value of the tuning parameter.

The adaptive algorithm studied in this paper (which corresponds to the process $\{Z_k\}$ defined below) requires for both its definition and study the introduction of an intermediate "stopped" process, which we now define:

First, in order to take into account potential jumps outside the space $\Theta$, we extend the parameter space with a cemetery point, $\theta_c \notin \Theta$, and define $\overline{\Theta} \stackrel{\text{def}}{=} \Theta \cup \{\theta_c\}$. It is convenient to introduce the family of transition kernels $\{Q_{\tilde{\gamma}}, \tilde{\gamma} \geq 0\}$ such that for any $\tilde{\gamma} \geq 0$, $(x, \theta) \in \mathsf{X} \times \Theta$, $A \in \mathcal{B}(\mathsf{X})$ and $B \in \mathcal{B}(\overline{\Theta})$,

$$
\begin{aligned}
Q_{\tilde{\gamma}}(x, \theta; A \times B) &= \int_A P_\theta(x, dy) \mathbb{1}\{\theta + \tilde{\gamma} H(\theta, y) \in B\} \\
&\quad + \delta_{\theta_c}(B) \int_A P_\theta(x, dy) \mathbb{1}\{\theta + \tilde{\gamma} H(\theta, y) \notin \Theta\},
\end{aligned}
\tag{10}
$$

where $\delta_\theta$ denotes the Dirac delta function at $\theta \in \Theta$. In its general form, the *basic* version of the adaptive MCMC algorithm considered here may be written as follows. Set $\theta_0 = \theta \in \Theta$, $X_0 = x \in \mathsf{X}$ and, for $k \geq 0$, define recursively the sequence $\{(X_k, \theta_k), k \geq 0\}$: if $\theta_k = \theta_c$, then set $\theta_{k+1} = \theta_c$ and $X_{k+1} = x$, otherwise $(X_{k+1}, \theta_{k+1}) \sim Q_{\rho_{k+1}}(X_k, \theta_k; \cdot)$, where $\boldsymbol{\rho} = \{\rho_k\}$ is a sequence of stepsizes. The sequence $\{(X_k, \theta_k)\}$ is a *nonhomogeneous* Markov chain on the product space $\mathsf{X} \times \overline{\Theta}$. This nonhomogeneous Markov chain defines a probability measure on the canonical state space $(\mathsf{X} \times \overline{\Theta})^{\mathbb{N}}$ equipped with the canonical product $\sigma$-algebra. We denote by $\mathcal{F} = \{\mathcal{F}_k, k \geq 0\}$ the natural filtration of this Markov chain and by $\mathbb{P}^{\boldsymbol{\rho}}_{x,\theta}$ and $\mathbb{E}^{\boldsymbol{\rho}}_{x,\theta}$ the probability and the expectation associated with this Markov chain starting from $(x, \theta) \in \mathsf{X} \times \Theta$.

Because of the interaction with feedback between $X_k$ and $\theta_k$, the stability of this inhomogeneous Markov chain is often difficult to establish. This is a long-standing problem in the field of stochastic optimization. Known practical remedies for this problem include the reprojections on a fixed set (see [23]) and the more recent reprojection on random varying boundaries proposed in [11, 12] and generalized in [3].

More precisely, we first define the notion of *compact coverage* of $\Theta$. A family of compact subsets $\{\mathcal{K}_q, q \geq 0\}$ of $\Theta$ is said to be a compact coverage if

(11) $$\bigcup_{q \geq 0} \mathcal{K}_q = \Theta \quad \text{and} \quad \mathcal{K}_q \subset \text{int}(\mathcal{K}_{q+1}), \qquad q \geq 0,$$

where $\text{int}(A)$ denotes the interior of set $A$. Let $\boldsymbol{\gamma} \stackrel{\text{def}}{=} \{\gamma_k\}$ be a monotone nonincreasing sequence of positive numbers and let $\mathsf{K}$ be a subset of $\mathsf{X}$. For a sequence $\mathbf{a} = \{a_k\}$ and an integer $l$, we define the "shifted" sequence $\mathbf{a}^{\leftarrow l}$ as follows: for any $k \geq 1$, $a_k^{\leftarrow l} \stackrel{\text{def}}{=} a_{k+l}$. Let $\Pi : \mathsf{X} \times \overline{\Theta} \to \mathsf{K} \times \mathcal{K}_0$ be a measurable function. Define the *homogeneous* Markov chain $Z_k = \{(X_k, \theta_k, \kappa_k, \nu_k), k \geq 0\}$ on the product space $\mathsf{Z} \stackrel{\text{def}}{=} \mathsf{X} \times \overline{\Theta} \times \mathbb{N} \times \mathbb{N}$, with transition probability $R : \mathsf{Z} \times \mathcal{B}(\mathsf{Z}) : \to [0, 1]$ algorithmically defined as follows (note that in order to alleviate notation, the dependence of $R$ on both $\boldsymbol{\gamma}$ and $\{\mathcal{K}_q, q \geq 0\}$ is implicit throughout the paper). For any $(x, \theta, \kappa, \nu) \in \mathsf{Z}$:

1. If $\nu = 0$, then draw $(X', \theta') \sim Q_{\gamma_\kappa}(\Pi(x, \theta); \cdot)$; otherwise, draw $(X', \theta') \sim Q_{\gamma_{\kappa+\nu}}(x, \theta; \cdot)$.
2. If $\theta' \in \mathcal{K}_\kappa$, then set $\kappa' = \kappa$ and $\nu' = \nu + 1$; otherwise, set $\kappa' = \kappa + 1$, and $\nu' = 0$.

In words, $\kappa$ and $\nu$ are counters: $\kappa$ is the index of the current active truncation set and $\nu$ counts the number of iterations since the last reinitialization. The event $\{\nu_k = 0\}$ indicates that a reinitialization occurs. The algorithm is restarted at iteration $k$ from a point in $\mathsf{K} \times \mathcal{K}_0$ with the "smaller" sequence of stepsizes $\boldsymbol{\gamma}^{\leftarrow \kappa_k}$. Note the important fact, at the heart of our analysis, that, between reinitializations, this process coincides with the *basic* version of the algorithm described earlier, with $\boldsymbol{\rho} = \boldsymbol{\gamma}^{\leftarrow \kappa_k}$. This is formalized in Lemma 1 below.

This algorithm is reminiscent of the projection on random varying boundaries proposed in [11, 12]: whenever the current iterate wanders outside the active truncation set, the algorithm is reinitialized with a smaller initial value of the stepsize and a larger truncation set.

The homogeneous Markov chain $\{Z_k, k \geq 0\}$ defines a probability measure on the canonical state space $\mathsf{Z}^\mathbb{N}$ equipped with the canonical product $\sigma$-algebra. We denote by $\mathcal{G} = \{\mathcal{G}_k, k \geq 0\}$, $\overline{\mathbb{P}}_{x,\theta,k,l}$ and $\overline{\mathbb{E}}_{x,\theta,k,l}$ the filtration, probability and expectation associated with this Markov chain started from $(x, \theta, k, l) \in \mathsf{Z}$. For simplicity, we will use the shorthand notation $\overline{\mathbb{E}}_\star$ and $\overline{\mathbb{P}}_\star$ for $\overline{\mathbb{E}}_{x,\theta} \stackrel{\text{def}}{=} \overline{\mathbb{E}}_{x,\theta,0,0}$ and $\overline{\mathbb{P}}_{x,\theta} \stackrel{\text{def}}{=} \overline{\mathbb{P}}_{x,\theta,0,0}$ for all $(x, \theta) \in \mathsf{X} \times \Theta$.

These probability measures depend upon the deterministic sequence $\boldsymbol{\gamma}$. The dependence will be implicit hereafter. We define recursively $\{T_n, n \geq 0\}$ the sequence of successive reinitialization times

(12) $$T_{n+1} = \inf\{k \geq T_n + 1, \nu_k = 0\} \qquad \text{with } T_0 = 0,$$

where, by convention, $\inf\{\varnothing\} = \infty$. It may be shown that, under mild conditions on $\{P_\theta, \theta \in \Theta\}$, $\{H(\theta, x), (\theta, x) \in \Theta \times \mathsf{X}\}$ and the sequence $\boldsymbol{\gamma}$,

$$\overline{\mathbb{P}}_\star\left(\sup_{n \geq 0} \kappa_n < \infty\right) = \overline{\mathbb{P}}_\star\left(\bigcup_{n=0}^{\infty} \{T_n = \infty\}\right) = 1,$$

that is, the number of reinitializations of the procedure described above is finite $\overline{\mathbb{P}}_\star$-a.s. We postpone to Sections 5, 6 and 7 the presentation and discussion of simple sufficient conditions that ensure that this holds in concrete situations. We will, however, assume this property to hold in Sections 3 and 4. Again, we stress the fact that our analysis of the homogeneous Markov chain $\{Z_k\}$ ("the algorithm") for a given sequence $\boldsymbol{\gamma}$ relies on the study of the inhomogeneous Markov chain defined earlier (the "stopped process"), for the sequences $\{\boldsymbol{\rho}_k \stackrel{\text{def}}{=} \boldsymbol{\gamma}^{\leftarrow \kappa_k}\}$ of stepsizes. It is therefore important to precisely and probabilistically relate these two processes. This is the aim of the lemma below (adapted from [3], Lemma 4.1).

Define, for $\mathcal{K} \subset \Theta$,

(13) $$\sigma(\mathcal{K}) = \inf\{k \geq 1, \theta_k \notin \mathcal{K}\}.$$

LEMMA 1. *Given any $m \geq 1$, any nonnegative measurable function $\Psi_m : (\mathsf{X} \times \overline{\Theta})^m \to \mathbb{R}^+$, for any integer $k \geq 0$ and $(x, \theta) \in \mathsf{X} \times \Theta$, satisfies*

$$\overline{\mathbb{E}}_{x,\theta,k,0}[\Psi_m(X_1, \theta_1, \ldots, X_m, \theta_m)\mathbb{1}\{T_1 \geq m\}]$$
$$= \mathbb{E}_{\Pi(x,\theta)}^{\boldsymbol{\gamma}^{\leftarrow k}}[\Psi_m(X_1, \theta_1, \ldots, X_m, \theta_m)\mathbb{1}\{\sigma(\mathcal{K}_k) \geq m\}].$$

**3. Law of large numbers.** Hereafter, for a probability distribution $\mathbb{P}$, the various kinds of convergence—in probability, almost-sure and weak (in distribution) are denoted, respectively, $\stackrel{\text{prob.}}{\longrightarrow}_{\mathbb{P}}$, $\stackrel{\text{a.s.}}{\longrightarrow}_{\mathbb{P}}$ and $\stackrel{\mathcal{D}}{\longrightarrow}_{\mathbb{P}}$.

3.1. *Assumptions.* As pointed out in the Introduction, an LLN has been obtained for a particular adaptive MCMC algorithm by Haario, Saksman and Tamminen [19], using mixingale theory (see [24]). Our approach is more in line with the martingale proof of the LLN for Markov chains and is based on the existence and regularity of the solutions of Poisson's equation and martingale limit theory. The existence and appropriate properties of those solutions can be easily established under a uniform (in the parameter $\theta$) version of the $V$-uniform ergodicity of the transition kernels $P_\theta$ [see condition (A1) below and Proposition 3].

We will use the following notation throughout the paper. For $W: \mathsf{X} \to [1, \infty)$ and $f: \mathsf{X} \to \mathbb{R}$ a measurable function, define

$$\|f\|_W = \sup_{x \in \mathsf{X}} \frac{|f(x)|}{W(x)} \quad \text{and} \quad \mathcal{L}_W = \{f : \|f\|_W < \infty\}. \tag{14}$$

We will also consider functions $f: \Theta \times \mathsf{X} \to \mathbb{R}$. We will often use the shorthand notation $f_\theta(x) = f(\theta, x)$ for all $\theta, x \in \Theta \times \mathsf{X}$ in order to avoid ambiguities. We will assume that $f_\theta \equiv 0$ whenever $\theta \notin \Theta$ except when $f_\theta$ does not depend on $\theta$, that is, $f_\theta \equiv f_{\theta'}$ for any $(\theta, \theta') \in \Theta \times \Theta$. Let $W: \mathsf{X} \to [1, \infty)$. We say that the family of functions $\{f_\theta : \mathsf{X} \to \mathbb{R}, \theta \in \Theta\}$ is $W$-*Lipschitz* if, for any compact subset $\mathcal{K} \subset \Theta$,

$$\sup_{\theta \in \mathcal{K}} \|f_\theta\|_W < \infty \quad \text{and} \quad \sup_{(\theta,\theta') \in \mathcal{K} \times \mathcal{K}, \theta \neq \theta'} |\theta - \theta'|^{-1} \|f_\theta - f_{\theta'}\|_W < \infty. \tag{15}$$

(A1) For any $\theta \in \Theta$, $P_\theta$ has $\pi$ as its stationary distribution. In addition, there exists a function $V: \mathsf{X} \to [1, \infty)$ such that $\sup_{x \in \mathsf{K}} V(x) < \infty$ (with $\mathsf{K}$ defined in Section 2) and such that, for any compact subset $\mathcal{K} \subset \Theta$:

(i) *Minorization condition.* There exist $\mathsf{C} \in \mathcal{B}(\mathsf{X})$, $\epsilon > 0$ and a probability measure $\varphi$ (all three depending on $\mathcal{K}$) such that $\varphi(\mathsf{C}) > 0$ and, for all $\mathsf{A} \in \mathcal{B}(\mathsf{X})$ and $\theta, x \in \mathcal{K} \times \mathsf{C}$,

$$P_\theta(x, \mathsf{A}) \geq \epsilon \varphi(\mathsf{A}).$$

(ii) *Drift condition.* There exist constants $\lambda \in [0, 1)$, $b \in (0, \infty)$ (depending on $V$, $\mathsf{C}$ and $\mathcal{K}$) satisfying

$$P_\theta V(x) \leq \begin{cases} \lambda V(x), & x \notin \mathsf{C}, \\ b, & x \in \mathsf{C}, \end{cases}$$

for all $\theta \in \mathcal{K}$.

(A2) For any compact subset $\mathcal{K} \subset \Theta$ and any $r \in [0, 1]$, there exists a constant $C$ (depending on $\mathcal{K}$ and $r$) such that, for any $(\theta, \theta') \in \mathcal{K} \times \mathcal{K}$ and $f \in \mathcal{L}_{V^r}$,

$$\|P_\theta f - P_{\theta'} f\|_{V^r} \leq C \|f\|_{V^r} |\theta - \theta'|,$$

where $V$ is given in (A1).

(A3) $\{H_\theta, \theta \in \Theta\}$ is $V^\beta$-Lipschitz for some $\beta \in [0, 1/2]$ with $V$ defined in (A1).

REMARK 1. Note that for the sake of clarity and simplicity, we restrict here our results to the case where $\{P_\theta, \theta \in \Theta\}$ satisfy one-step drift and minorization conditions. As shown in [3], the more general case where either an $m$-step drift or minorization condition is assumed for $m > 1$ requires one to modify the algorithm in order to prevent large jumps in the parameter space (see [2, 3]). This mainly leads to substantial notational complications, but the arguments remain essentially unchanged.

Conditions of type (A1) to establish geometric ergodicity have been extensively studied over the last decade for the Metropolis–Hastings algorithms. Typically, the required drift function depends on the target distribution $\pi$, which makes our requirement of uniformity in $\theta \in \mathcal{K}$ in (A1) reasonable and relatively easy to establish (see Sections 6 and 7). The following theorem, due to [29] and recently improved by [8], converts information about the drift and minorization condition into information about the long-term behavior of the chain:

THEOREM 2. *Assume* (A1). *Then, for all* $\theta \in \Theta$, $P_\theta$ *admits* $\pi$ *as its unique stationary probability measure, and* $\pi(V) < \infty$. *Let* $\mathcal{K} \subset \Theta$ *be a compact subset and* $r \in [0,1]$. *There exists* $\tilde{\rho} < 1$ *depending only (and explicitly) on the constants* $r$, $\epsilon$, $\varphi(\mathsf{C})$, $\lambda$ *and* $b$ [*given in* (A1)] *such that, whenever* $\rho \in (\tilde{\rho}, 1)$, *there exists* $C < \infty$ *depending only (and explicitly) on* $r$, $\rho$, $\epsilon$, $\varphi(\mathsf{C})$ *and* $b$ *such that for any* $f \in \mathcal{L}_{V^r}$, *all* $\theta \in \mathcal{K}$ *and* $k \geq 0$,

$$\|P_\theta^k f - \pi(f)\|_{V^r} \leq C \|f\|_{V^r} \rho^k. \tag{16}$$

Formulas for $\rho$ and $C$ are given in [29], Theorem 2.3, and have been later improved in [8], Section 2.1.

This theorem automatically ensures the existence of solutions to Poisson's equation. More precisely, for all $\theta, x \in \Theta \times \mathsf{X}$ and $f \in \mathcal{L}_{V^r}$, $\sum_{k=0}^\infty |P_\theta^k f(x) - \pi(f)| < \infty$ and $\hat{f}_\theta \stackrel{\text{def}}{=} \sum_{k=0}^\infty (P_\theta^k f - \pi(f))$ is a solution of Poisson's equation

$$\hat{f}_\theta - P_\theta \hat{f}_\theta = f - \pi(f). \tag{17}$$

Poisson's equation has proven to be a fundamental tool for the analysis of additive functionals, in particular for establishing limit theorems such as the (functional) central limit theorem (see, e.g., [9, 13, 18, 30], [28], Chapter 17.

The Lipschitz continuity of the transition kernel $P_\theta$ as a function of $\theta$ [assumption (A2)] does not seem to have been studied for the Metropolis–Hastings algorithm. We establish this continuity for the SRWM algorithm and the independent MH algorithm (IMH) in Sections 6 and 7. This assumption, used in conjunction with (A1), allows one to establish the Lipschitz continuity of the solution of Poisson's equation.

PROPOSITION 3. *Assume* (A1). *Suppose that the family of functions* $\{f_\theta, \theta \in \Theta\}$ *is* $V^r$-*Lipschitz, for some* $r \in [0,1]$. *Define, for any* $\theta \in \Theta$, $\hat{f}_\theta \stackrel{\text{def}}{=} \sum_{k=0}^\infty (P_\theta^k f_\theta - \pi(f_\theta))$. *Then, for any compact set* $\mathcal{K}$, *there exists a constant* $C$ *such that, for any* $(\theta, \theta') \in \mathcal{K}$,

$$\|\hat{f}_\theta\|_{V^r} + \|P_\theta \hat{f}_\theta\|_{V^r} \leq C \sup_{\theta \in \mathcal{K}} \|f_\theta\|_{V^r}, \tag{18}$$

$$\|\hat{f}_\theta - \hat{f}_{\theta'}\|_{V^r} + \|P_\theta \hat{f}_\theta - P_{\theta'} \hat{f}_{\theta'}\|_{V^r} \leq C |\theta - \theta'| \sup_{\theta \in \mathcal{K}} \|f_\theta\|_{V^r}. \tag{19}$$

The proof is given in Appendix B.

REMARK 2. The regularity of the solutions of Poisson's equation has been studied, under various ergodicity and regularity conditions on the mapping $\theta \mapsto P_\theta$ (see, e.g., [9] and [7] for regularity under conditions implying $V$-uniform geometric ergodicity). The results of the proposition above are sharper than those reported in the literature, because all the transition kernels $P_\theta$ share the same limiting distribution $\pi$, a property which plays a key role in the proof.

We finish this section with a convergence result for the chain $\{X_k\}$ under the probability $\mathbb{P}^{\boldsymbol{\rho}}_{x,\theta}$, which is an important and direct byproduct of the property mentioned in the remark immediately above. This result improves on [19, 21] and [6].

PROPOSITION 4. *Assume* (A1)–(A3), *let* $\rho \in (0,1)$ *be as in* (16), *let* $\bar{\boldsymbol{\rho}} = \{\bar{\rho}_k\}$ *be a positive and finite, nonincreasing sequence such that* $\limsup_{k \to \infty} \bar{\rho}_{k-n(k)}/\bar{\rho}_k < +\infty$ *with, when* $k \geq k_0$ *for some integer* $k_0$, $n(k) := \lfloor \log(\bar{\rho}_k/\bar{\rho})/\log(\rho) \rfloor + 1 \leq k$, *and* $n(k) = 0$ *otherwise, where* $\bar{\rho} \in [\rho_1, +\infty)$. *Let* $f \in \mathcal{L}_{V^{1-\beta}}$, *where* $V$ *is defined in* (A1), *and* $\beta$ *is defined in* (A3). *Let* $\mathcal{K}$ *be a compact subset of* $\Theta$. *Then, there exists a constant* $C \in (0, +\infty)$ [*depending only on* $\mathcal{K}$, *the constants in* (A1) *and* $\bar{\boldsymbol{\rho}}$] *such that, for any* $(x, \theta) \in \mathsf{X} \times \mathcal{K}$,

$$|\mathbb{E}^{\bar{\boldsymbol{\rho}}}_{x,\theta}\{(f(X_k) - \pi(f))\mathbb{1}\{\sigma(\mathcal{K}) \geq k\}\}| \leq C\|f\|_{V^{1-\beta}}\bar{\rho}_k V(x).$$

The proof appears in Appendix B.

3.2. *Law of large numbers.* We prove in this section a law of large numbers (LLN) under $\overline{\mathbb{P}}_\star$ for $n^{-1}\sum_{k=1}^n f_{\theta_k}(X_k)$, where $\{f_\theta, \theta \in \Theta\}$ is a set of sufficiently regular functions. It is worth noting here that we need not require that the sequence $\{\theta_k\}$ converges in order to establish our result. The proof is based on the identity

$$f_{\theta_k}(X_k) - \int_\mathsf{X} \pi(dx) f_{\theta_k}(x) = \hat{f}_{\theta_k}(X_k) - P_{\theta_k}\hat{f}_{\theta_k}(X_k),$$

where $\hat{f}_\theta$ is a solution of Poisson's equation (17). The decomposition

$$\begin{aligned}
&\hat{f}_{\theta_k}(X_k) - P_{\theta_k}\hat{f}_{\theta_k}(X_k) \\
(20) \quad &= (\hat{f}_{\theta_{k-1}}(X_k) - P_{\theta_{k-1}}\hat{f}_{\theta_{k-1}}(X_{k-1})) \\
&\quad + (\hat{f}_{\theta_k}(X_k) - \hat{f}_{\theta_{k-1}}(X_k)) + (P_{\theta_{k-1}}\hat{f}_{\theta_{k-1}}(X_{k-1}) - P_{\theta_k}\hat{f}_{\theta_k}(X_k)),
\end{aligned}$$

displays the different terms that need to be controlled to prove the LLN. The first term in the decomposition is (except at the time of a jump) a martingale difference sequence. As we shall see, this is the leading term in the decomposition, and the other terms are remainders which are easily dealt with, thanks to the regularity of the solutions of Poisson's equation under (A1). The term $\hat{f}_{\theta_k}(X_k) - \hat{f}_{\theta_{k-1}}(X_k)$ can be interpreted as the perturbation introduced by adaptation. We preface our main result, Theorem 8, with two intermediate propositions concerned with the control of the fluctuations of the sum $\sum_{k=1}^n (f_{\theta_k}(X_k) - \int_X \pi(dx) f_{\theta_k}(x))$ for the inhomogeneous chain $\{(X_k, \theta_k)\}$ under the probability $\mathbb{P}^{\boldsymbol{\rho}}_{x,\theta}$. The following lemma, whose proof is given in Appendix A, is required in order to prove Propositions 6 and 7:

LEMMA 5. *Assume* (A1). *Let* $\mathcal{K} \subset \Theta$ *be a compact set and* $r \in [0,1]$ *be a constant. There exists a constant* $C$ [*depending only on* $r$, $\mathcal{K}$ *and the constants in* (A1)] *such that, for any sequences* $\boldsymbol{\rho} = \{\rho_k\}$ *and* $\mathbf{a} = \{a_k\}$ *of positive numbers and for any* $(x, \theta) \in \mathsf{X} \times \mathcal{K}$,

$$\mathbb{E}^{\boldsymbol{\rho}}_{x,\theta}[V^r(X_k) \mathbb{1}\{\sigma(\mathcal{K}) \geq k\}] \leq CV^r(x), \tag{21}$$

$$\mathbb{E}^{\boldsymbol{\rho}}_{x,\theta}\left[\max_{1 \leq m \leq k}(a_m V(X_m))^r \mathbb{1}\{\sigma(\mathcal{K}) \geq m\}\right] \leq C^r \left(\sum_{m=1}^k a_m\right)^r V^r(x), \tag{22}$$

$$\mathbb{E}^{\boldsymbol{\rho}}_{x,\theta}\left[\max_{1 \leq m \leq n} \mathbb{1}\{\sigma(\mathcal{K}) \geq m\} \sum_{k=1}^m V^r(X_k)\right] \leq CnV^r(x). \tag{23}$$

PROPOSITION 6. *Assume* (A1)–(A3). *Let* $\{f_\theta, \theta \in \Theta\}$ *be a* $V^\alpha$-*Lipschitz family of functions for some* $\alpha \in [0, 1-\beta)$, *where* $V$ *is defined in* (A1) *and* $\beta$ *is defined in* (A3). *Let* $\mathcal{K}$ *be a compact subset of* $\Theta$. *Then, for any* $p \in (1, 1/(\alpha + \beta)]$, *there exists a constant* $C$ [*depending only on* $p$, $\mathcal{K}$ *and the constants in* (A1)] *such that, for any sequence* $\boldsymbol{\rho} = \{\rho_k\}$ *of positive numbers satisfying* $\sum_{k=1}^\infty k^{-1} \rho_k < \infty$, *we have, for all* $(x, \theta) \in \mathsf{X} \times \mathcal{K}$, $\delta > 0$, $\tilde{\gamma} > \alpha$ *and integer* $l \geq 1$,

$$\mathbb{P}^{\boldsymbol{\rho}}_{x,\theta}\left[\sup_{m \geq l}(m^{-1}|M_m|) \geq \delta\right]$$
$$\leq C\delta^{-p} \sup_{\theta \in \mathcal{K}} \|f_\theta\|^p_{V^\alpha} l^{-\{(p/2) \wedge (p-1)\}} V^{\alpha p}(x), \tag{24}$$

$$\mathbb{P}^{\boldsymbol{\rho}}_{x,\theta}\left[\sup_{m \geq l}\left\{\mathbb{1}\{\sigma(\mathcal{K}) > m\} m^{-\tilde{\gamma}} \left|\sum_{k=1}^m \left(f_{\theta_k}(X_k) - \int_\mathsf{X} f_{\theta_k}(x)\pi(dx)\right) - M_m\right|\right\} \geq \delta\right]$$
$$\leq C\delta^{-p} \sup_{\theta \in \mathcal{K}} \|f_\theta\|^p_{V^\alpha} \left\{\left(\sum_{k=1}^\infty (k \vee l)^{-\tilde{\gamma}} \rho_k\right)^p V(x)^{\beta p} + l^{-(\tilde{\gamma} - \alpha)p}\right\} \tag{25}$$

$$\times V^{\alpha p}(x),$$

where $\sigma(\mathcal{K})$ is given in (13), $\hat{f}_\theta \stackrel{\text{def}}{=} \sum_{k=1}^\infty [P_\theta^k f_\theta - \pi(f_\theta)]$ is a solution of Poisson's equation (17) and

$$(26) \qquad M_m \stackrel{\text{def}}{=} \mathbb{1}\{\sigma(\mathcal{K}) > m\} \sum_{k=1}^m [\hat{f}_{\theta_{k-1}}(X_k) - P_{\theta_{k-1}} \hat{f}_{\theta_{k-1}}(X_{k-1})].$$

REMARK 3. The result provides us with some useful insights into the properties of MCMC with vanishing adaptation. First, whenever $\{\theta_k\} \subset \mathcal{K} \subset \Theta$ for a deterministic compact set $\mathcal{K}$, the bounds above give explicit rates of convergence for ergodic averages (1). The price we must pay for adaptation is apparent on the righ-hand side of (25), as reported in [1]. The constraints on $p, \beta$ and $\alpha$ illustrate the tradeoff between the rate of convergence, the smoothness of adaptation and the class of functions covered by our result. Scenarios of interest include the case where assumption (A3) is satisfied with $\beta = 0$ [in other words, for any compact subset $\mathcal{K} \subset \Theta$, the function $\sup_{\theta \in \mathcal{K}} \sup_{x \in \mathsf{X}} |H_\theta(x)| < \infty$] or the case where (A3) holds for any $\beta \in (0, 1/2]$, both of which imply that the results of the proposition hold for any $\alpha < 1$ (see Theorem 15 for an application).

PROOF OF PROPOSITION 6. For notational simplicity, we set $\sigma \stackrel{\text{def}}{=} \sigma(\mathcal{K})$. Let $p \in (1, 1/(\alpha + \beta)]$ and $\mathcal{K} \subset \Theta$ be a compact set. In this proof, $C$ is a constant which only depends on the constants in (A1)–(A3), $p$ and $\mathcal{K}$; this constant may take different values upon each appearance. Theorem 2 implies that for any $\theta \in \Theta$, $\hat{f}_\theta$ exists and is a solution of Poisson's equation (17). We decompose the sum $\mathbb{1}\{\sigma > m\} \sum_{k=1}^m (f_{\theta_k}(X_k) - \int_{\mathsf{X}} f_{\theta_k}(x) \pi(dx))$ as $M_m + R_m^{(1)} + R_m^{(2)}$, where

$$R_m^{(1)} \stackrel{\text{def}}{=} \mathbb{1}\{\sigma > m\} \sum_{k=1}^m (\hat{f}_{\theta_k}(X_k) - \hat{f}_{\theta_{k-1}}(X_k)),$$

$$R_m^{(2)} \stackrel{\text{def}}{=} \mathbb{1}\{\sigma > m\} (P_{\theta_0} \hat{f}_{\theta_0}(X_0) - P_{\theta_m} \hat{f}_{\theta_m}(X_m)).$$

We consider these terms separately. First, since $\mathbb{1}\{\sigma > m\} = \mathbb{1}\{\sigma > m\} \mathbb{1}\{\sigma \geq k\}$ for $0 \leq k \leq m$,

$$|M_m| = \mathbb{1}\{\sigma > m\} \left| \sum_{k=1}^m (\hat{f}_{\theta_{k-1}}(X_k) - P_{\theta_{k-1}} \hat{f}_{\theta_{k-1}}(X_{k-1})) \mathbb{1}\{\sigma \geq k\} \right| \leq |\widetilde{M}_m|,$$

where

$$\widetilde{M}_m \stackrel{\text{def}}{=} \sum_{k=1}^m [\hat{f}_{\theta_{k-1}}(X_k) - P_{\theta_{k-1}} \hat{f}_{\theta_{k-1}}(X_{k-1})] \mathbb{1}\{\sigma \geq k\}.$$

By Proposition 3, equation (18), there exists a constant $C$ such that, for all $\theta \in \mathcal{K}$, $\|\hat{f}_\theta\|_{V^\alpha} + \|P_\theta \hat{f}_\theta\|_{V^\alpha} \leq C \sup_{\theta \in \mathcal{K}} \|f_\theta\|_{V^\alpha}$. Since $0 \leq \alpha p \leq 1$, by using (21) in Lemma 5, we have, for all $x, \theta \in \mathsf{X} \times \mathcal{K}$,

$$
\begin{aligned}
(27) \quad & \mathbb{E}^{\boldsymbol{\rho}}_{x,\theta}\{(|\hat{f}_{\theta_{k-1}}(X_k)|^p + |P_{\theta_{k-1}}\hat{f}_{\theta_{k-1}}(X_{k-1})|^p)\mathbb{1}\{\sigma \geq k\}\} \\
& \leq CV^{\alpha p}(x)\sup_{\theta \in \mathcal{K}}\|f_\theta\|^p_{V^\alpha}.
\end{aligned}
$$

Since

$$
\begin{aligned}
& \mathbb{E}^{\boldsymbol{\rho}}_{x,\theta}\{[\hat{f}_{\theta_{k-1}}(X_k) - P_{\theta_{k-1}}\hat{f}_{\theta_{k-1}}(X_{k-1})]\mathbb{1}\{\sigma \geq k\}|\mathcal{F}_{k-1}\} \\
& \quad = (P_{\theta_{k-1}}\hat{f}_{\theta_{k-1}}(X_{k-1}) - P_{\theta_{k-1}}\hat{f}_{\theta_{k-1}}(X_{k-1}))\mathbb{1}\{\sigma \geq k\} = 0,
\end{aligned}
$$

$\{\widetilde{M}_m\}$ is a $(\mathbb{P}^{\boldsymbol{\rho}}_{x,\theta}, \{\mathcal{F}_k\})$-adapted martingale with increments bounded in $L^p$. Using Burkholder's inequality for $p > 1$ ([20], Theorem 2.10), we have

$$
\begin{aligned}
(28) \quad & \mathbb{E}^{\boldsymbol{\rho}}_{x,\theta}\{|\widetilde{M}_m|^p\} \\
& \leq C_p \mathbb{E}^{\boldsymbol{\rho}}_{x,\theta}\bigg\{\bigg(\sum_{k=1}^m |\hat{f}_{\theta_{k-1}}(X_k) - P_{\theta_{k-1}}\hat{f}_{\theta_{k-1}}(X_{k-1})|^2 \mathbb{1}\{\sigma \geq k\}\bigg)^{p/2}\bigg\},
\end{aligned}
$$

where $C_p$ is a universal constant. For $p \geq 2$, by Minkowski's inequality,

$$
\begin{aligned}
(29) \quad & \mathbb{E}^{\boldsymbol{\rho}}_{x,\theta}\{|\widetilde{M}_m|^p\} \\
& \leq C_p\bigg\{\sum_{k=1}^m (\mathbb{E}^{\boldsymbol{\rho}}_{x,\theta}\{|\hat{f}_{\theta_{k-1}}(X_k) - P_{\theta_{k-1}}\hat{f}_{\theta_{k-1}}(X_{k-1})|^p \mathbb{1}\{\sigma \geq k\}\})^{2/p}\bigg\}^{p/2}.
\end{aligned}
$$

For $1 \leq p \leq 2$, we have

$$
\begin{aligned}
(30) \quad & \mathbb{E}^{\boldsymbol{\rho}}_{x,\theta}\bigg\{\bigg(\sum_{k=1}^m |\hat{f}_{\theta_{k-1}}(X_k) - P_{\theta_{k-1}}\hat{f}_{\theta_{k-1}}(X_{k-1})|^2 \mathbb{1}\{\sigma \geq k\}\bigg)^{p/2}\bigg\} \\
& \leq \mathbb{E}^{\boldsymbol{\rho}}_{x,\theta}\bigg\{\sum_{k=1}^m |\hat{f}_{\theta_{k-1}}(X_k) - P_{\theta_{k-1}}\hat{f}_{\theta_{k-1}}(X_{k-1})|^p \mathbb{1}\{\sigma \geq k\}\bigg\}.
\end{aligned}
$$

By combining the two cases above and using (27), we obtain for any $x, \theta \in \mathsf{X} \times \mathcal{K}$ and $p > 1$,

$$
(31) \quad \mathbb{E}^{\boldsymbol{\rho}}_{x,\theta}\{|\widetilde{M}_m|^p\} \leq Cm^{(p/2) \vee 1} \sup_{\theta \in \mathcal{K}}\|f_\theta\|^p_{V^\alpha} V^{\alpha p}(x).
$$

Let $l \geq 1$. By Birnbaum and Marshall's [10] inequality (a straightforward adaptation of Birnbaum and Marshall [10] result is given in Proposition 22,

Appendix A) and (31), there exists a constant $C$ such that

$$\mathbb{P}_{x,\theta}^{\boldsymbol{\rho}}\left[\sup_{m\geq l}(m^{-1}|\widetilde{M}_m|)\geq \delta\right] \leq \delta^{-p}\left\{\sum_{m=l}^{\infty}(m^{-p}-(m+1)^{-p})\mathbb{E}_{x,\theta}^{\boldsymbol{\rho}}\{|\widetilde{M}_m|^p\}\right\}$$

$$\leq C\delta^{-p}\sup_{\theta\in\mathcal{K}}\|f_\theta\|_{V^\alpha}^p\left(\sum_{m=l}^{\infty}m^{-p-1+(p/2)\vee 1}\right)V^{\alpha p}(x)$$

$$\leq C\delta^{-p}\sup_{\theta\in\mathcal{K}}\|f_\theta\|_{V^\alpha}^p l^{-\{(p/2)\wedge(p-1)\}}V^{\alpha p}(x),$$

which proves (24). We now consider $R_m^{(1)}$. Equation (19) shows that there exists a constant $C$ such that, for any $(\theta,\theta')\in\mathcal{K}\times\mathcal{K}$, $\|\hat{f}_\theta - \hat{f}_{\theta'}\|_{V^\alpha} \leq C|\theta - \theta'|\sup_{\theta\in\mathcal{K}}\|f_\theta\|_{V^\alpha}$. On the other hand, by construction, $\theta_k - \theta_{k-1} = \rho_k H(\theta_{k-1}, X_k)$ for $\sigma \geq k$ and, under assumption (A3), there exists a constant $C$ such that, for any $x,\theta \in \mathsf{X}\times\mathcal{K}$, $|H(\theta,x)| \leq CV^\beta(x)$. Therefore, there exists $C$ such that, for any $m \geq l$ and $\tilde{\gamma} > \alpha$,

$$m^{-\tilde{\gamma}}|R_m^{(1)}| = m^{-\tilde{\gamma}}\mathbb{1}\{\sigma > m\}\left|\sum_{k=1}^{m}\{\hat{f}_{\theta_k}(X_k) - \hat{f}_{\theta_{k-1}}(X_k)\}\right|$$

$$\leq C\sup_{\theta\in\mathcal{K}}\|f_\theta\|_{V^\alpha}\sum_{k=1}^{m}(k\vee l)^{-\tilde{\gamma}}\rho_k V^{\alpha+\beta}(X_k)\mathbb{1}\{\sigma\geq k\}.$$

Hence, using Minkowski's inequality and (21), one deduces that there exists $C$ such that, for $(x,\theta) \in \mathsf{X}\times\mathcal{K}$, $l \geq 1$ and $\tilde{\gamma} > \alpha$,

$$(32) \quad \begin{aligned}&\mathbb{E}_{x,\theta}^{\boldsymbol{\rho}}\left\{\sup_{m\geq l}m^{-\tilde{\gamma}p}|R_m^{(1)}|^p\right\}\\ &\leq C\sup_{\theta\in\mathcal{K}}\|f_\theta\|_{V^\alpha}^p\left(\sum_{k=1}^{\infty}(k\vee l)^{-\tilde{\gamma}}\rho_k\right)^p V^{(\alpha+\beta)p}(x).\end{aligned}$$

Consider now $R_m^{(2)}$. The term $P_{\theta_0}\hat{f}_{\theta_0}(X_0)$ does not pose any problem. From Lemma 5 [equation (22)] there exists a constant $C$ such that, for all $x,\theta \in \mathsf{X}\times\mathcal{K}$ and $0 < \alpha < \tilde{\gamma}$,

$$(33) \quad \begin{aligned}&\mathbb{E}_{x,\theta}^{\boldsymbol{\rho}}\left\{\sup_{m\geq l}m^{-\tilde{\gamma}p}|P_{\theta_m}\hat{f}_{\theta_m}(X_m)|^p\mathbb{1}\{\sigma > m\}\right\}\\ &\leq C\sup_{\theta\in\mathcal{K}}\|f_\theta\|_{V^\alpha}^p\mathbb{E}_{x,\theta}^{\boldsymbol{\rho}}\left\{\sup_{m\geq l}|m^{-\tilde{\gamma}/\alpha}V^\alpha(X_m)|^p\mathbb{1}\{\sigma>m\}\right\}\\ &\leq C\sup_{\theta\in\mathcal{K}}\|f_\theta\|_{V^\alpha}^p\left(\sum_{k=l}^{\infty}k^{-\tilde{\gamma}/\alpha}\right)^{\alpha p}V^{\alpha p}(x).\end{aligned}$$

The case $\alpha = 0$ is straightforward. From Markov's inequality, (32) and (33), one deduces (25). $\square$

We can now apply the results of Proposition 6 for the inhomogeneous Markov chain defined below (10) to the time-homogeneous time chain $\{Z_k\}$ under the assumption that the number of reinitializations $\kappa_n$ is $\overline{\mathbb{P}}_\star$ almost surely finite. Note that the very general form of the result will allow us to prove a central limit theorem in Section 4.

PROPOSITION 7. *Let $\{\mathcal{K}_q, q \geq 0\}$ be a compact coverage of $\Theta$ and let $\gamma = \{\gamma_k\}$ be a nonincreasing positive sequence such that $\sum_{k=1}^\infty k^{-1}\gamma_k < \infty$. Consider the time-homogeneous Markov chain $\{Z_k\}$ on $\mathsf{Z}$ with transition probability $R$, as defined in Section 2. Assume* (A1)–(A3), *and let $\mathcal{F} \stackrel{\mathrm{def}}{=} \{f_\theta, \theta \in \Theta\}$ be a $V^\alpha$-Lipschitz family of functions for some $\alpha \in [0, 1-\beta)$, with $\beta$ as in* (A3) *and $V$ as in* (A1). *Assume, in addition, that $\overline{\mathbb{P}}_\star\{\lim_{n\to\infty} \kappa_n < \infty\} = 1$. Then, for any $f_\theta \in \mathcal{F}$,*

$$(34) \qquad n^{-1}\sum_{k=1}^n \left(f_{\theta_k}(X_k) - \int_\mathsf{X} f_{\theta_k}(x)\pi(dx)\right)\mathbb{1}\{\nu_k \neq 0\} \xrightarrow{\mathrm{a.s.}}_{\overline{\mathbb{P}}_\star} 0.$$

PROOF. Without loss of generality, we may assume that, for any $\theta \in \Theta$, $\int_\mathsf{X} f_\theta(x)\pi(dx) = 0$. Let $p \in (1, 1/(\alpha+\beta)]$ and define $\kappa_\infty \stackrel{\mathrm{def}}{=} \lim_{n\to\infty} \kappa_n$. By construction of the $T_k$'s [see (12)], and since $\kappa_\infty < \infty$ $\overline{\mathbb{P}}_\star$-a.s., we deduce that $T_{\kappa_\infty} < \infty$, $\overline{\mathbb{P}}_\star$-a.s. We now decompose $S_n = \sum_{k=1}^n f_{\theta_k}(X_k)\mathbb{1}\{\nu_k \neq 0\}$ as $S_n \stackrel{\mathrm{def}}{=} S_n^{(1)} + S_n^{(2)}$, where $S_n^{(1)} \stackrel{\mathrm{def}}{=} \sum_{k=1}^{T_{\kappa_\infty} \wedge n} f_{\theta_k}(X_k)\mathbb{1}\{\nu_k \neq 0\}$ and $S_n^{(2)} \stackrel{\mathrm{def}}{=} \sum_{k=T_{\kappa_\infty}+1}^n f_{\theta_k}(X_k)$. Since $T_{\kappa_\infty} < \infty$ $\overline{\mathbb{P}}_\star$-a.s., $\sum_{k=1}^{T_{\kappa_\infty}} |f_{\theta_k}(X_k)|\mathbb{1}\{\nu_k \neq 0\} < \infty$, $\overline{\mathbb{P}}_\star$-a.s., showing that $n^{-1}S_n^{(1)} \to 0$, $\overline{\mathbb{P}}_\star$-a.s. We now bound the second term. For any integers $n$ and $K$,

$$(35) \quad \begin{aligned} &\overline{\mathbb{P}}_\star\left[\sup_{m \geq n} m^{-1}|S_m^{(2)}| \geq \delta, \kappa_\infty \leq K\right] \\ &\leq \sum_{i=0}^K \overline{\mathbb{P}}_\star\left[\sup_{m\geq n} m^{-1}\left|\sum_{k=T_i+1}^m f_{\theta_k}(X_k)\right| \geq \delta, \kappa_\infty = i, T_i \leq n/2\right] \\ &\quad + \overline{\mathbb{P}}_\star[T_{\kappa_\infty} > n/2]. \end{aligned}$$

Since, for $\kappa_\infty = i$, $T_{i+1} = \infty$,

$$\left\{\sup_{m\geq n} m^{-1}\left|\sum_{k=T_i+1}^m f_{\theta_k}(X_k)\right| \geq \delta, \kappa_\infty = i, T_i \leq n/2\right\}$$

$$\subset \left\{\sup_{m\geq n} \mathbb{1}\{T_{i+1} > m\}m^{-1}\left|\sum_{k=T_i+1}^m f_{\theta_k}(X_k)\right| \geq \delta, T_i \leq n/2\right\}$$

$$\subset \left\{ \left( \sup_{m \geq n/2} \mathbb{1}\{\sigma(\mathcal{K}_i) > m\} m^{-1} \left| \sum_{k=1}^{m} f_{\theta_k}(X_k) \right| \right) \circ \tau^{T_i} \geq \delta \right\},$$

where we have used the fact that $T_{i+1} = T_i + \sigma(\mathcal{K}_i) \circ \tau^{T_i}$, where $\tau$ is the shift operator on the canonical space of the chain $\{Z_n\}$. As a consequence, by applying Lemma 1 [noting that $\mathbb{1}\{\sigma(\mathcal{K}_i) \geq m\}\mathbb{1}\{\sigma(\mathcal{K}_i) > m\} = \mathbb{1}\{\sigma(\mathcal{K}_i) > m\}$ and $\mathbb{1}\{\sigma(\mathcal{K}_i) > m\} \in \mathcal{F}_m$], the strong Markov property, Proposition 6 with $\tilde{\gamma} = 1$, the fact that $\{\gamma_k\}$ is nonincreasing and the fact that $\sup_{x \in \mathsf{K}} V(x) < \infty$, we have, for $0 \leq i \leq K$,

$$\overline{\mathbb{P}}_\star \left[ \sup_{m \geq n} m^{-1} \left| \sum_{k=T_i+1}^{m} f_{\theta_k}(X_k) \right| \geq \delta, \kappa_\infty = i, T_i \leq n/2 \bigg| \mathcal{G}_{T_i} \right]$$

(36)
$$\leq \mathbb{P}_{\Pi(X_{T_i}, \theta_{T_i})}^{\gamma^{\leftarrow i}} \left[ \sup_{m \geq n/2} \mathbb{1}\{\sigma(\mathcal{K}_i) > m\} m^{-1} \left| \sum_{k=1}^{m} f_{\theta_k}(X_k) \right| \geq \delta \right]$$

$$\leq C \delta^{-p} \sup_{\theta \in \mathcal{K}_i} \|f_\theta\|_{V^\alpha}^p \left\{ \left( \sum_{k=1}^{\infty} [\lfloor n/2 \rfloor \vee k]^{-1} \gamma_k \right)^p + \left( \frac{n}{2} \right)^{-p(1-\alpha)} + \left( \frac{n}{2} \right)^{-\{(p/2) \wedge (p-1)\}} \right\}.$$

By Kronecker's lemma, the condition $\sum_{k=1}^{\infty} k^{-1} \gamma_k < \infty$ implies that $n^{-1} \sum_{k=1}^{n} \gamma_k \to 0$ as $n \to \infty$, showing that

$$\sum_{k=1}^{\infty} [\lfloor n/2 \rfloor \vee k]^{-1} \gamma_k \leq (\lfloor n/2 \rfloor)^{-1} \sum_{k=1}^{\lfloor n/2 \rfloor} \gamma_k + \sum_{k=\lfloor n/2 \rfloor+1}^{\infty} k^{-1} \gamma_k \to 0,$$

as $n \to \infty$. Combining this with (35) and (36) shows that, for any $K, \delta, \eta > 0$, there exists $N$ such that, for $n \geq N$,

$$\overline{\mathbb{P}}_\star \left[ \sup_{m \geq n} m^{-1} |S_m^{(2)}| \geq \delta, \kappa_\infty \leq K \right] \leq \eta.$$

Now, for $K$ large enough that $\overline{\mathbb{P}}_\star[\kappa_\infty > K] \leq \eta$, the result above shows that there exists an $N$ such that, for any $n \geq N$, $\overline{\mathbb{P}}_\star[\sup_{m \geq n} m^{-1} |S_m^{(2)}| \geq \delta] \leq 2\eta$, concluding the proof. $\square$

REMARK 4. How one checks $\overline{\mathbb{P}}_\star(\lim_{n \to \infty} \kappa_n < \infty) = 1$ depends on the particular algorithm used to update the parameters. Verifiable conditions have been established in [3] for checking the stability of the algorithm; see Sections 5, 6 and 7.

We may now state our main consistency result.

THEOREM 8. *Let $\{\mathcal{K}_q, q \geq 0\}$ be a compact coverage of $\Theta$ and let $\gamma = \{\gamma_k\}$ be a nonincreasing positive sequence such that $\sum_{k=1}^{\infty} k^{-1} \gamma_k < \infty$. Consider the time-homogeneous Markov chain $\{Z_k\}$ on $\mathsf{Z}$ with transition probability $R$, as defined in Section 2. Assume (A1)–(A3) and let $f:\mathsf{X} \to \mathbb{R}$ be a function such that $\|f\|_{V^\alpha} < \infty$ for some $\alpha \in [0, 1-\beta)$, with $\beta$ as in (A3) and $V$ as in (A1). Assume, in addition, that $\overline{\mathbb{P}}_\star\{\lim_{n\to\infty} \kappa_n < \infty\} = 1$. Then,*

$$(37) \qquad n^{-1} \sum_{k=1}^{n} [f(X_k) - \pi(f)] \xrightarrow{\text{a.s.}}_{\overline{\mathbb{P}}_\star} 0.$$

PROOF. We may assume that $\pi(f) = 0$. From Proposition 7, it is sufficient to prove that $n^{-1} \sum_{j=1}^{\kappa_n} |f(X_{T_j})| \xrightarrow{\text{a.s.}}_{\overline{\mathbb{P}}_\star} 0$. Since $\kappa_\infty < \infty$ $\overline{\mathbb{P}}_\star$-a.s., $\sum_{j=1}^{\kappa_\infty} |f(X_{T_j})| < \infty$ $\overline{\mathbb{P}}_\star$-a.s. The proof follows from $\sum_{j=1}^{\kappa_n} |f(X_{T_j})| \leq \sum_{j=1}^{\kappa_\infty} |f(X_{T_j})|$. □

**4. Invariance principle.** We shall now prove an invariance principle. As in the case of homogeneous Markov chains, more stringent conditions are required here than for the simple LLN. In particular, we will require here that the series $\{\theta_k\}$ converges $\overline{\mathbb{P}}_\star$-a.s. This is in contrast with simple consistency for which boundedness of $\{\theta_k\}$ was sufficient. The main idea of the proof consists of approximating $n^{-1/2} \sum_{k=1}^{n} \{f(X_k) - \pi(f)\}$ with a triangular array of martingale differences sequence, and then applying an invariance principle for martingale differences to show the desired result.

THEOREM 9. *Let $\{\mathcal{K}_q, q \geq 0\}$ be a compact coverage of $\Theta$ and let $\gamma = \{\gamma_k\}$ be a nonincreasing positive sequence such that $\sum_{k=1}^{\infty} k^{-1/2} \gamma_k < \infty$. Consider the time-homogeneous Markov chain $\{Z_k\}$ on $\mathsf{Z}$ with transition probability $R$, as defined in Section 2. Assume (A1)–(A3) and let $f:\mathsf{X} \to \mathbb{R}$ satisfy $f \in \mathcal{L}_{V^\alpha}$, where $V$ is defined in (A3) and $\alpha \in [0, (1-\beta)/2)$ with $\beta$ as in (A1). Define, for any $\theta \in \Theta$,*

$$(38) \qquad \sigma^2(\theta, f) \stackrel{\text{def}}{=} \pi[(\hat{f}_\theta - P_\theta \hat{f}_\theta)^2] \qquad \text{with } \hat{f}_\theta \stackrel{\text{def}}{=} \sum_{k=0}^{\infty} [P_\theta^k f - \pi(f)].$$

*Assume, in addition, that there exists a random variable $\theta_\infty \in \Theta$, such that $\int_\mathsf{X} \pi(dx) \hat{f}_{\theta_\infty}^2(x) < \infty$ and $\int_\mathsf{X} \pi(dx) (P_{\theta_\infty} \hat{f}_{\theta_\infty}(x))^2 < \infty$ $\overline{\mathbb{P}}_\star$-a.s. and*

$$\limsup_{n\to\infty} |\theta_n - \theta_\infty| = 0, \qquad \overline{\mathbb{P}}_\star\text{-a.s.}$$

*Then,*

$$n^{-1/2} \sum_{k=1}^{n} [f(X_k) - \pi(f)] \xrightarrow{\mathcal{D}}_{\overline{\mathbb{P}}_\star} Z,$$

where the random variable $Z$ has characteristic function $\overline{\mathbb{E}}_\star[\exp(-\frac{1}{2}\sigma^2(\theta_\infty, f)t^2)]$. If in addition $\sigma(\theta_\infty, f) > 0$, $\overline{\mathbb{P}}_\star$-a.s., then

$$(39) \quad \frac{1}{\sqrt{n}\sigma(\theta_\infty, f)} \sum_{k=1}^{n}(f(X_k) - \pi(f)) \xrightarrow{\mathcal{D}}_{\overline{\mathbb{P}}_\star} \mathcal{N}(0,1).$$

PROOF. Without loss of generality, we suppose that $\pi(f) = 0$. The proof again relies on a martingale approximation. Set, for $k \geq 1$,

$$(40) \quad \xi_k \stackrel{\text{def}}{=} [\hat{f}_{\theta_{k-1}}(X_k) - P_{\theta_{k-1}}\hat{f}_{\theta_{k-1}}(X_{k-1})]\mathbb{1}\{\nu_{k-1} \neq 0\}.$$

Since $f$ is $V^\alpha$-Lipschitz, Proposition 3 shows that $\{\hat{f}_\theta, \theta \in \Theta\}$ and $\{P_\theta \hat{f}_\theta, \theta \in \Theta\}$ are $V^\alpha$-Lipschitz. Since $2\alpha < 1$, this implies that $\{P_\theta \hat{f}_\theta^2, \theta \in \Theta\}$ and $\{(P_\theta \hat{f}_\theta)^2, \theta \in \Theta\}$ are $V^{2\alpha}$-Lipschitz. We deduce that $\{\xi_k\}$ is a $(\overline{\mathbb{P}}_\star, \{\mathcal{G}_k, k \geq 0\})$-adapted square-integrable martingale difference sequence, that is, for all $k \geq 1$, $\overline{\mathbb{E}}_\star[\xi_k^2] < \infty$ and $\overline{\mathbb{E}}_\star[\xi_k|\mathcal{G}_{k-1}] = 0$, $\overline{\mathbb{P}}_\star$-a.s. We are going to prove that with $Z$ a r.v. with characteristic function $\overline{\mathbb{E}}_\star[\exp(-\frac{1}{2}\sigma^2(\theta_\infty, f)t^2)]$,

$$(41) \quad \frac{1}{\sqrt{n}} \sum_{k=1}^{n} \xi_k \xrightarrow{\mathcal{D}}_{\overline{\mathbb{P}}_\star} Z,$$

$$(42) \quad \frac{1}{\sqrt{n}} \sum_{k=1}^{n} f(X_k) - \frac{1}{\sqrt{n}} \sum_{k=1}^{n} \xi_k \xrightarrow{\text{a.s.}}_{\overline{\mathbb{P}}_\star} 0.$$

To show (41), we use [20], Corollary 3.1 of Theorem 3.2. We need to establish that:

(a) the sequence $n^{-1}\sum_{k=1}^{n} \overline{\mathbb{E}}_\star[\xi_k^2|\mathcal{G}_{k-1}]$ converges in $\overline{\mathbb{P}}_\star$-probability to $\sigma^2(\theta_\infty, f)$;

(b) the conditional Lindeberg condition is satisfied, that is,

$$(43) \quad \text{for all } \varepsilon > 0 \quad n^{-1}\sum_{k=1}^{n} \overline{\mathbb{E}}_\star[\xi_k^2 \mathbb{1}\{|\xi_k| \geq \varepsilon\sqrt{n}\}|\mathcal{G}_{k-1}] \xrightarrow{\text{prob.}}_{\overline{\mathbb{P}}_\star} 0.$$

We first prove (a). Note that

$$\overline{\mathbb{E}}_\star[\xi_k^2|\mathcal{G}_{k-1}] = [P_{\theta_{k-1}}\hat{f}_{\theta_{k-1}}^2(X_{k-1}) - (P_{\theta_{k-1}}\hat{f}_{\theta_{k-1}}(X_{k-1}))^2]\mathbb{1}\{\nu_{k-1} \neq 0\}.$$

Since $\{P_\theta \hat{f}_\theta^2, \theta \in \Theta\}$ and $\{(P_\theta \hat{f}_\theta)^2, \theta \in \Theta\}$ are $V^{2\alpha}$-Lipschitz and $2\alpha \in [0, 1 - \beta)$, we may apply Proposition 7 to prove that

$$\frac{1}{n}\sum_{k=1}^{n}[P_{\theta_{k-1}}\hat{f}_{\theta_{k-1}}^2(X_{k-1}) - (P_{\theta_{k-1}}\hat{f}_{\theta_{k-1}}(X_{k-1}))^2]\mathbb{1}\{\nu_{k-1} \neq 0\}$$

$$- \frac{1}{n}\sum_{k=0}^{n-1}\int_{\mathsf{X}} \pi(dx)[P_{\theta_k}\hat{f}_{\theta_k}^2(x) - (P_{\theta_k}\hat{f}_{\theta_k}(x))^2]\mathbb{1}\{\nu_k \neq 0\} \xrightarrow{\text{a.s.}}_{\overline{\mathbb{P}}_\star} 0.$$

For any $j \geq 0$ and $\kappa_\infty = j$, $\{\theta_k, k > T_j\} \subset \mathcal{K}_j$, which, together with the $V^{2\alpha}$-Lipschitz property and the dominated convergence theorem, implies that

$$\lim_{k \to \infty} \mathbb{1}\{\kappa_\infty = j\} \int_\mathsf{X} \pi(dx) \left[ P_{\theta_k} \hat{f}^2_{\theta_k}(x) - (P_{\theta_k} \hat{f}_{\theta_k}(x))^2 \right]$$
$$= \mathbb{1}\{\kappa_\infty = j\} \int_\mathsf{X} \pi(dx) [P_{\theta_\infty} \hat{f}^2_{\theta_\infty}(x) - (P_{\theta_\infty} \hat{f}_{\theta_\infty}(x))^2], \qquad \overline{\mathbb{P}}_\star\text{-a.s.}$$

By the dominated convergence theorem and the fact that $\overline{\mathbb{P}}_\star(\kappa_\infty < \infty) = 1$,

$$\lim_{k \to \infty} \int_\mathsf{X} \pi(dx) [P_{\theta_k} \hat{f}^2_{\theta_k}(x) - (P_{\theta_k} \hat{f}_{\theta_k}(x))^2] \mathbb{1}\{\nu_k \neq 0\}$$
$$= \lim_{k \to \infty} \sum_{j=0}^\infty \left\{ \int_\mathsf{X} \pi(dx) [P_{\theta_k} \hat{f}^2_{\theta_k}(x) - (P_{\theta_k} \hat{f}_{\theta_k}(x))^2] \right\} \mathbb{1}\{\nu_k \neq 0\} \mathbb{1}\{\kappa_\infty = j\}$$
$$= \sum_{j=0}^\infty \left\{ \lim_{k \to \infty} \int_\mathsf{X} \pi(dx) [P_{\theta_k} \hat{f}^2_{\theta_k}(x) - (P_{\theta_k} \hat{f}_{\theta_k}(x))^2] \right\} \mathbb{1}\{\kappa_\infty = j\}, \qquad \overline{\mathbb{P}}_\star\text{-a.s.,}$$

and the Cesàro convergence theorem finally shows that

$$\lim_{n \to \infty} \frac{1}{n} \sum_{k=0}^{n-1} \int_\mathsf{X} \pi(dx) [P_{\theta_k} \hat{f}^2_{\theta_k}(x) - (P_{\theta_k} \hat{f}_{\theta_k}(x))^2] \mathbb{1}\{\nu_k \neq 0\}$$
$$= \int_\mathsf{X} \pi(dx) [P_{\theta_\infty} \hat{f}^2_{\theta_\infty}(x) - (P_{\theta_\infty} \hat{f}_{\theta_\infty}(x))^2], \qquad \overline{\mathbb{P}}_\star\text{-a.s.}$$

We now establish the conditional Lindeberg condition in (b). We use the following lemma, which is a conditional version of [14], Lemma 3.3.

LEMMA 10. *Let $\mathcal{G}$ be a $\sigma$-field and $X$ a random variable such that $\mathbb{E}[X^2|\mathcal{G}] < \infty$. Then, for any $\varepsilon > 0$,*

$$4\mathbb{E}[|X|^2 \mathbb{1}\{|X| \geq \varepsilon\} | \mathcal{G}] \geq \mathbb{E}[|X - \mathbb{E}[X|\mathcal{G}]|^2 \mathbb{1}\{|X - \mathbb{E}[X|\mathcal{G}]| \geq 2\varepsilon\} | \mathcal{G}].$$

Using Dvoretzky's lemma, we have for any $\varepsilon, M > 0$ and $n$ sufficiently large,

$$n^{-1} \sum_{k=1}^n \overline{\mathbb{E}}_\star[\xi_k^2 \mathbb{1}\{|\xi_k| \geq \varepsilon\sqrt{n}\} | \mathcal{G}_{k-1}]$$
$$\leq 4n^{-1} \sum_{k=0}^{n-1} \int_\mathsf{X} P_{\theta_k}(X_k, dx) \hat{f}^2_{\theta_k}(x) \mathbb{1}\{|\hat{f}_{\theta_k}(x)| \geq M\} \mathbb{1}\{\nu_k \neq 0\}, \qquad \overline{\mathbb{P}}_\star\text{-a.s.}$$

Proceeding as above, the right-hand side of the previous display converges $\overline{\mathbb{P}}_\star$-a.s. to

$$\int_\mathsf{X} \pi(dx) \hat{f}^2_{\theta_\infty}(x) \mathbb{1}\{|\hat{f}_{\theta_\infty}(x)| \geq M\},$$

where we have used the fact that, for any $\theta \in \Theta$, $\pi P_\theta = \pi$. Since $\int_X \pi(dx) \times \hat{f}^2_{\theta_\infty}(x) < \infty$ $\overline{\mathbb{P}}_\star$-a.s., the monotone convergence theorem implies that

$$\lim_{M \to \infty} \int_X \pi(dx) \hat{f}^2_{\theta_\infty}(x) \mathbb{1}\{|\hat{f}_{\theta_\infty}(x)| \geq M\} = 0, \qquad \overline{\mathbb{P}}_\star\text{-a.s.},$$

showing that the conditional Lindeberg condition (b) holds.

In order to prove equation (42), we proceed along the lines of the proof of Proposition 7. First, since $\kappa_\infty < \infty$, $\overline{\mathbb{P}}_\star$-a.s., $\sum_{k=1}^{T_{\kappa_\infty}} |f(X_k)| + \sum_{k=1}^{T_{\kappa_\infty}} |\xi_k| < \infty$, $\overline{\mathbb{P}}_\star$-a.s., which implies that

$$(44) \qquad n^{-1/2} \left| \sum_{k=1}^{T_{\kappa_\infty}} f(X_k) - \sum_{k=1}^{T_{\kappa_\infty}} \xi_k \right| \xrightarrow{\text{a.s.}}_{\overline{\mathbb{P}}_\star} 0.$$

Second, proceeding as in the proof of (36) and using (25) with $\tilde{\gamma} = 1/2$ and some $p \in (1, 1/(\alpha + \beta)]$, since $f$ is $V^\alpha$-Lipschitz, we have that for any $0 \leq i \leq K$ for some $K > 0$ and $n > 0$,

$$(45) \quad \begin{aligned}&\overline{\mathbb{P}}_\star\left[\sup_{m \geq n} m^{-1/2} \left| \sum_{k=T_i+1}^{m} f(X_k) - \sum_{k=T_i+1}^{m} \xi_k \right| \geq \delta, \kappa_\infty = i, T_i \leq n/2 \bigg| \mathcal{G}_{T_i} \right] \\ &\leq C\delta^{-p} \|f\|_{V^\alpha}^p \left\{ \left( \sum_{k=1}^{\infty} [\lfloor n/2 \rfloor \vee k]^{-1/2} \gamma_k \right)^p + \left(\frac{n}{2}\right)^{-p(1/2-\alpha)} \right\}.\end{aligned}$$

Under the assumption $\sum_{k=1}^{\infty} k^{-1/2} \gamma_k < \infty$, proceeding as below equation (36), one can show that $\lim_{n \to \infty} \sum_{k=1}^{\infty} [\lfloor n/2 \rfloor) \vee k]^{-1/2} \gamma_k = 0$. Arguing as in (35), we conclude that

$$(46) \qquad m^{-1/2} \left| \sum_{k=T_i+1}^{m} f(X_k) - \sum_{k=T_i+1}^{m} \xi_k \right| \xrightarrow{\text{a.s.}}_{\overline{\mathbb{P}}_\star} 0.$$

The proof of (42) follows from (44) and (46). The proof of (39) follows from [20], Corollary 3.2 of Theorem 3.3. □

**5. Stability and convergence of the stochastic approximation process.** In order to conclude the part of this paper dedicated to the general theory of adaptive MCMC algorithms, we now present generally verifiable conditions under which the number of reinitializations of the algorithm that produces the Markov chain $\{Z_k\}$ described in Section 2 is $\overline{\mathbb{P}}_\star$-a.e. finite. This is a difficult problem *per se*, which has been worked out in a companion paper, [3]. We here briefly introduce the conditions under which this key property is satisfied and give (without proof) the main stability result. The reader should refer to [3] for more details.

As mentioned in the Introduction, the convergence of the stochastic approximation procedure is closely related to the stability of the noiseless sequence $\bar{\theta}_{k+1} = \bar{\theta}_k + \gamma_{k+1} h(\bar{\theta}_k)$. A practical technique for proving the stability of the noiseless sequence consists, when possible, of determining a Lyapunov function $w : \Theta \to [0, \infty)$ such that $\langle \nabla w(\theta), h(\theta) \rangle \leq 0$, where $\nabla w$ denotes the gradient of $w$ with respect to $\theta$ and, for $u, v \in \mathbb{R}^n$, $\langle u, v \rangle$ is their Euclidean inner product (we will later on also use the notation $|v| = \sqrt{\langle v, v \rangle}$ to denote the Euclidean norm of $v$). This indeed shows that the noiseless sequence $\{w(\bar{\theta}_k)\}$ eventually decreases, showing that $\lim_{k \to \infty} w(\bar{\theta}_k)$ exists. It should therefore not be surprising if such a Lyapunov function can play an important role in showing the stability of the noisy sequence $\{\theta_k\}$. With this in mind, we can now detail the conditions required to prove our convergence result:

(A4) $\Theta$ is an open subset of $\mathbb{R}^{n_\theta}$. The mean field $h : \Theta \to \mathbb{R}^{n_\theta}$ is continuous, and there exists a continuously differentiable function $w : \Theta \to [0, \infty)$ [with the convention $w(\theta) = \infty$ when $\theta \notin \Theta$] such that:

(i) For any $M > 0$, the level set $\mathcal{W}_M \stackrel{\text{def}}{=} \{\theta \in \Theta, w(\theta) \leq M\} \subset \Theta$ is compact;

(ii) the set of stationary point(s) $\mathcal{L} \stackrel{\text{def}}{=} \{\theta \in \Theta, \langle \nabla w(\theta), h(\theta) \rangle = 0\}$ belongs to the interior of $\Theta$;

(iii) for any $\theta \in \Theta$, $\langle \nabla w(\theta), h(\theta) \rangle \leq 0$ and the closure of $w(\mathcal{L})$ has an empty interior.

Finally we require some conditions on the sequence of stepsizes $\boldsymbol{\gamma} = \{\gamma_k\}$.

(A5) The sequence $\boldsymbol{\gamma} = \{\gamma_k\}$ is nonincreasing, positive and

$$\sum_{k=1}^{\infty} \gamma_k = \infty \quad \text{and} \quad \sum_{k=1}^{\infty} \{\gamma_k^2 + k^{-1/2} \gamma_k\} < \infty.$$

The following theorem is a straightforward simplification of [3], Theorems 5.4 and 5.5, and shows that the tail probability of the number of reinitializations decreases faster than any exponential, and that the parameter sequence $\{\theta_k\}$ converges to the stationary set $\mathcal{L}$. For a point $x$ and a set $A$ we define $d(x, A) \stackrel{\text{def}}{=} \inf\{|x - y| : y \in A\}$.

THEOREM 11. *Let $\{\mathcal{K}_q, q \geq 0\}$ be a compact coverage of $\Theta$ and let $\boldsymbol{\gamma} = \{\gamma_k\}$ be a real-valued sequence. Consider the time-homogeneous Markov chain $\{Z_k\}$ on $\mathsf{Z}$ with transition probability $R$, as defined in Section 2. Assume (A1)–(A5). Then,*

$$\limsup_{k \to \infty} k^{-1} \log \left( \sup_{(x, \theta) \in \mathsf{X} \times \Theta} \overline{\mathbb{P}}_{x, \theta} \left[ \sup_{n \geq 0} \kappa_n \geq k \right] \right) = -\infty,$$

$$\inf_{(x,\theta)\in\mathsf{X}\times\Theta} \overline{\mathbb{P}}_{x,\theta}\left[\lim_{k\to\infty} d(\theta_k,\mathcal{L})=0\right]=1.$$

**6. Consistency and invariance principle for the adaptive N-SRW kernel.** In this section we show how our results can be applied to the adaptive N-SRWM algorithm proposed by Haario, Saksman and Tamminen [19] and described in Section 1. We first illustrate how the conditions required to prove the LLN in [19] can be alleviated. In particular, no boundedness condition is required on the parameter set $\Theta$, but rather conditions on the tails of the target distribution $\pi$. We then extend these results further and prove a central limit theorem (Theorem 15).

In view of the results proved above it is required:

(a) to prove the ergodicity and regularity conditions for the Markov kernels outlined in assumption (A1);

(b) to prove that the reinitializations occur finitely many times (stability) and that $\{\theta_k\}$ eventually converges. Note again that the convergence property is only required for the CLT.

We first focus on (a). The geometric ergodicity of the SRWM kernel has been studied by Roberts and Tweedie [31] and refined by Jarner and Hansen [22]; the regularity of the SRWM kernel has not, to the best of our knowledge, been considered in the literature. The geometric ergodicity of the SRWM kernel mainly depends on the tail properties of the target distribution $\pi$. We will therefore restrict our discussion to target distributions that satisfy the following set of conditions. These are not minimal, but easy to check in practice (see [22] for details).

(M) The probability density $\pi$ is defined on $\mathsf{X}=\mathbb{R}^{n_x}$ for some integer $n_x$ and has the following properties:

  (i) It is bounded, bounded away from zero on every compact set and continuously differentiable.

  (ii) It is super-exponential, that is,

$$\lim_{|x|\to+\infty}\left\langle \frac{x}{|x|},\nabla\log\pi(x)\right\rangle=-\infty.$$

  (iii) The contours $\partial\mathsf{A}(x)=\{y:\pi(y)=\pi(x)\}$ are asymptotically regular, that is,

$$\lim_{|x|\to+\infty}\sup\left\langle \frac{x}{|x|},\frac{\nabla\pi(x)}{|\nabla\pi(x)|}\right\rangle<0.$$

We now establish uniform minorization and drift conditions for the SRWM algorithm defined in (3). Let $\mathcal{M}(\mathsf{X})$ denote the set of probability densities

w.r.t. the Lebesgue measure $\lambda^{\text{Leb}}$. For any $a, b > 0$, define $\mathcal{Q}_{a,b}(\mathsf{X}) \subset \mathcal{M}(\mathsf{X})$ as follows:

$$\text{(47)} \qquad \mathcal{Q}_{a,b}(\mathsf{X}) \stackrel{\text{def}}{=} \left\{ q \in \mathcal{M}(\mathsf{X}), q(x) = q(-x) \text{ and } \inf_{|x| \leq a} q(x) \geq b \right\}.$$

PROPOSITION 12. *Assume* (M). *For any* $\eta \in (0,1)$, *set* $V = \pi^{-\eta}/(\sup_{x \in \mathsf{X}} \pi(x))^{-\eta}$. *Then:*

1. *For any nonempty compact set* $\mathsf{C} \subset \mathsf{X}$, *there exists* $a > 0$ *such that, for any* $b > 0$ *such that* $\mathcal{Q}_{a,b}(\mathsf{X}) \neq \varnothing$, *there exists* $\epsilon > 0$ *such that* $\mathsf{C}$ *is a* $(1, \epsilon)$-*small set for the elements of* $\{P_q^{\text{SRW}} : q \in \mathcal{Q}_{a,b}(\mathsf{X})\}$, *with minorization probability distribution* $\varphi$ *such that, for any* $\mathsf{A} \in \mathcal{B}(\mathsf{X})$, $\varphi(\mathsf{A}) = \lambda^{\text{Leb}}(\mathsf{A} \cap \mathsf{C})/\lambda^{\text{Leb}}(\mathsf{C})$, *that is,*

$$\text{(48)} \qquad \inf_{q \in \mathcal{Q}_{a,b}(\mathsf{X})} P_q^{\text{SRW}}(x, \mathsf{A}) \geq \epsilon \varphi(\mathsf{A}) \qquad \text{for all } x \in \mathsf{C} \text{ and } \mathsf{A} \in \mathcal{B}(\mathsf{X}).$$

2. *Furthermore, for any* $a > 0$ *and* $b > 0$ *such that* $\mathcal{Q}_{a,b}(\mathsf{X}) \neq \varnothing$,

$$\text{(49)} \qquad \sup_{q \in \mathcal{Q}_{a,b}(\mathsf{X})} \limsup_{|x| \to +\infty} \frac{P_q^{\text{SRW}} V(x)}{V(x)} < 1,$$

$$\text{(50)} \qquad \sup_{(x,q) \in \mathsf{X} \times \mathcal{Q}_{a,b}(\mathsf{X})} \frac{P_q^{\text{SRW}} V(x)}{V(x)} < +\infty.$$

3. *Let* $q, q' \in \mathcal{M}(\mathsf{X})$ *be two symmetric probability distributions. Then, for any* $r \in [0, 1]$ *and any* $f \in \mathcal{L}_{V^r}$, *we have*

$$\text{(51)} \qquad \|P_q^{\text{SRW}} f - P_{q'}^{\text{SRW}} f\|_{V^r} \leq 2\|f\|_{V^r} \int_{\mathsf{X}} |q(x) - q'(x)| \lambda^{\text{Leb}}(dx).$$

The proof appears in Appendix C.

As an example of an application, one can again consider the adaptive N-SRWM introduced earlier in Section 1, where the proposal distribution is $\mathcal{N}(0, \Gamma)$. In the following lemma, we show that the mapping $\Gamma \to P_{\mathcal{N}(0,\Gamma)}^{\text{SRW}}$ is Lipschitz continuous. This result can be generalized to distributions in the curved exponential family (see Proposition 16).

LEMMA 13. *Let* $\mathcal{K}$ *be a convex compact subset of* $\mathcal{C}_+^{n_x}$ *and set* $V = \pi^{-\eta}/(\sup_{\mathsf{X}} \pi)^{-\eta}$ *for some* $\eta \in (0,1)$. *For any* $r \in [0,1]$, *any* $\Gamma, \Gamma' \in \mathcal{K} \times \mathcal{K}$ *and any* $f \in \mathcal{L}_{V^r}$, *we have*

$$\|P_{\mathcal{N}(0,\Gamma)}^{\text{SRW}} f - P_{\mathcal{N}(0,\Gamma')}^{\text{SRW}} f\|_{V^r} \leq \frac{2n_x}{\lambda_{\min}(\mathcal{K})} \|f\|_{V^r} |\Gamma - \Gamma'|,$$

*where, for* $\Gamma \in \mathcal{C}_+^{n_x}$, $|\Gamma|^2 = \text{Tr}[\Gamma \Gamma^{\mathsf{T}}]$ *and* $\lambda_{\min}(\mathcal{K})$ *is the minimum possible eigenvalue for matrices in* $\mathcal{K}$.

The proof appears in Appendix D. We now turn to proving that the stochastic approximation procedure outlined by Haario, Saksman and Tamminen [19] is ultimately pathwise bounded and eventually converges. In the case of the algorithm proposed by Haario, Saksman and Tamminen [19], the parameter estimates $\mu_k$ and $\Gamma_k$ take the form of maximum likelihood estimates under the i.i.d. multivariate Gaussian model. It therefore comes as no surprise if the Lyapunov function $w$ required to check (A4) is the Kullback–Leibler divergence between the target density $\pi$ and the normal density $\mathcal{N}(\mu, \Gamma)$,

$$(52) \qquad w(\mu, \Gamma) = \log \det \Gamma + (\mu - \mu_\pi)^{\mathsf{T}} \Gamma^{-1} (\mu - \mu_\pi) + \mathrm{Tr}(\Gamma^{-1} \Gamma_\pi),$$

where $\mu_\pi$ and $\Gamma_\pi$ are the mean and covariance of the target distribution, defined in (9). Using straightforward algebra and the definition (8) of the mean field $h$, one can check that

$$\begin{aligned}(53)\quad &\langle \nabla w(\mu, \Gamma), h(\mu, \Gamma) \rangle \\ &= -2(\mu - \mu_\pi)^{\mathsf{T}} \Gamma^{-1} (\mu - \mu_\pi) \\ &\quad - \mathrm{Tr}(\Gamma^{-1}(\Gamma - \Gamma_\pi)\Gamma^{-1}(\Gamma - \Gamma_\pi)) - ((\mu - \mu_\pi)^{\mathsf{T}} \Gamma^{-1} (\mu - \mu_\pi))^2,\end{aligned}$$

that is, $\langle \nabla w(\theta), h(\theta) \rangle \leq 0$ for any $\theta \stackrel{\text{def}}{=} (\mu, \Gamma) \in \Theta$, with equality if and only if $\Gamma = \Gamma_\pi$ and $\mu = \mu_\pi$. The situation in this case is simple, as the set of stationary points $\{\theta \in \Theta, h(\theta) = 0\}$ is reduced to a single point, and the Lyapunov function $w$ goes to infinity as $|\mu| \to \infty$ or $\Gamma$ goes to the boundary of the cone of positive matrices.

It can now be shown that these results lead to the following intermediate lemma; see [3] for details.

LEMMA 14. *Assume* (M) *and let $H$ and $h$ be as in* (5) *and* (8). *Then,* (A3) *and* (A4) *are satisfied with $V = \pi^{-1}/(\sup_{\mathsf{X}} \pi)^{-1}$ for any $\beta \in (0, 1/2]$ and $w$ as in* (52). *In addition, the set of stationary points $\mathcal{L} \stackrel{\text{def}}{=} \{\theta \in \Theta \stackrel{\text{def}}{=} \mathbb{R}^{n_x} \times \mathcal{C}_+^{n_x}, \langle \nabla w(\theta), h(\theta) \rangle = 0\}$ is reduced to a single point $\theta_\pi = (\mu_\pi, \Gamma_\pi)$, whose components are respectively the mean and covariance of the distribution $\pi$.*

From Proposition 12 and Lemma 13, we deduce our main theorem for this section, concerned with the adaptive N-SRWM of [19] as described in Section 1, but with reprojections as in Section 2.

THEOREM 15. *Consider the process $\{Z_k\}$ with $\{P_\theta, \theta = (\mu, \Gamma) \in \Theta \stackrel{\text{def}}{=} \mathbb{R}^{n_x} \times \mathcal{C}_+^{n_x}, q_\theta = \mathcal{N}(0, \lambda \Gamma), \lambda > 0\}$ as in* (3), *$\{H_\theta, \theta \in \Theta\}$ as in equation* (5), *$\pi$ satisfying* (M), *$\boldsymbol{\gamma} = \{\gamma_k, k \geq 0\}$ satisfying* (A5) *and $\mathsf{K}$ a compact set. Let $W \stackrel{\text{def}}{=} \pi^{-1}/(\sup \pi)^{-1}$. Then, for any $\alpha \in [0, 1)$:*

1. *For any $f \in \mathcal{L}(W^\alpha)$ a strong LLN holds, that is,*

$$n^{-1} \sum_{k=1}^{n} \left( f(X_k) - \int_{\mathsf{X}} f(x)\pi(dx) \right) \xrightarrow{\text{a.s.}}_{\overline{\mathbb{P}}_\star} 0. \tag{54}$$

2. *For any $f \in \mathcal{L}(W^{\alpha/2})$ a CLT holds, that is,*

$$\frac{1}{\sqrt{n}} \sum_{k=1}^{n} [f(X_k) - \pi(f)] \xrightarrow{\mathcal{D}}_{\overline{\mathbb{P}}_\star} \mathcal{N}(0, \sigma^2(\theta_\pi, f)),$$

*if $\sigma(\theta_\pi, f) > 0$,*

$$\frac{1}{\sqrt{n}\sigma(\theta_\pi, f)} \sum_{k=1}^{n} (f(X_k) - \pi(f)) \xrightarrow{\mathcal{D}}_{\overline{\mathbb{P}}_\star} \mathcal{N}(0, 1),$$

*where $\theta_\pi = (\mu_\pi, \Gamma_\pi)$ and $\sigma^2(\theta_\pi, f)$ are defined in* (38).

The proof is immediate. We refer the reader to [19] for applications of this type of algorithm to various settings.

## 7. Application: matching $\pi$ with mixtures.

7.1. *Setup.* The independence Metropolis–Hastings algorithm (IMH) corresponds to the case where the proposal distribution used in an MH transition probability does not depend on the current state of the MCMC chain, that is, $q(x, y) = q(y)$ for some density $q \in \mathcal{M}(\mathsf{X})$. The transition kernel of the Metropolis algorithm is then given for $x \in \mathsf{X}$ and $A \in \mathcal{B}(\mathsf{X})$ by

$$P_q^{\text{IMH}}(x, A) = \int_A \alpha_q(x, y) q(y) \lambda^{\text{Leb}}(dy)$$
$$+ \mathbb{1}_A(x) \int_{\mathsf{X}} (1 - \alpha_q(x, y)) q(y) \lambda^{\text{Leb}}(dy) \tag{55}$$

$$\text{with } \alpha_q(x, y) = 1 \wedge \frac{\pi(y)q(x)}{\pi(x)q(y)}.$$

Irreducibility of Markov chains built on this model naturally require that $q(x) > 0$ whenever $\pi(x) > 0$. In fact, the performance of the IMH is known to depend on how well the proposal distribution mimics the target distribution, and this can be quantified in several ways. For example, it has been shown in [26] that the IMH sampler is geometrically ergodic if and only if there exists $\varepsilon > 0$ such that $q \in \mathcal{Q}_{\varepsilon,\pi} \subset \mathcal{M}(\mathsf{X})$, where

$$\mathcal{Q}_{\varepsilon,\pi} = \{q \in \mathcal{M}(\mathsf{X}) : \lambda^{\text{Leb}}(\{x \in \mathsf{X} : q(x)/\pi(x) < \varepsilon\}) = 0\}. \tag{56}$$

This condition implies that the whole state space $\mathsf{X}$ is a $(1, \varepsilon)$-small set, which in turn implies that convergence occurs uniformly, at a geometric rate bounded above by $1 - \varepsilon$. Given a family of candidate proposal distributions $\{q_\theta \in \mathcal{M}(\mathsf{X}), \theta \in \Theta\}$, it therefore seems natural to maximise $\theta \to \inf_{x \in \mathsf{X}} \pi(x)/q_\theta(x)$. However, although theoretically attractive, the optimization of this uniform criterion might be a very ambitious task in practice. Furthermore, it might not necessarily be a good choice *for a given parametric family of proposal distributions*: one might in this case try to optimize the transition probability for pathological features of $\pi$ with small probability under $\pi$, at the expense of more fundamental characteristics of the target, such as its global shape. Additionally, such pathological features can very often be taken care of by other specialized MCMC updates. Instead of this uniform criterion, we suggest the optimization of an average property of the ratio $\pi(x)/q_\theta(x)$ under $\pi$, which possesses the advantage of being more amenable to computation. It is argued in [15] that minimizing the total variation distance $\|\pi - q_\theta\|_{\mathrm{TV}}$ is a sensible criterion to optimize, since it can be proved that the expected acceptance probability is bounded below by $1 - \|\pi - q_\theta\|_{\mathrm{TV}}$, and that, for a bounded function $f$, the first covariance coefficient of the Markov chain in the stationary regime is bounded as follows: $\mathrm{cov}_\pi(f(X_k), f(X_{k+1})) \leq (5/2)^2 \sup_{x \in \mathsf{X}} |f| \|\pi - q_\theta\|_{\mathrm{TV}}$. However, no systematic way of effectively minimizing this criterion is described. We propose here to use the Kullback–Leibler divergence between the target distribution $\pi$ and an auxiliary distribution $\tilde{q}_\theta$ close in some sense to $q_\theta$,

$$(57) \qquad K(\pi\|\tilde{q}_\theta) = \int_\mathsf{X} \pi(x) \log \frac{\pi(x)}{\tilde{q}_\theta(x)} \lambda^{\mathrm{Leb}}(dx).$$

The proposal distribution $q_\theta$ of the IMH algorithm is then constructed from $\tilde{q}_\theta$. As we shall see, this offers an additional degree of freedom which, in particular, will be a simple way of ensuring that $\{q_\theta, \theta \in \Theta\} \subset \mathcal{Q}_{\varepsilon,\pi}$, defined in (56), for some $\varepsilon > 0$ (see Remark 7). The use of this criterion possesses several advantages. First, invoking Pinsker's inequality, it is possible to repeat Gasemyr [15] arguments. Second, it formalizes several ideas that have been proposed in the literature (cf. [15] and [17] among others). In [17] it is suggested to use the EM (Expectation–Minimization) algorithm in order to fit a mixture of normals in the, possibly penalized, maximum likelihood sense to samples from a preliminary run of an MCMC algorithm. This mixture can then be used to define the proposal distribution of an IMH. As we shall see, the choice of the Kullback–Leibler (KL) divergence corresponds precisely to this choice and naturally leads to an on-line EM algorithm that allows us to adjust $q_\theta$ to $\pi$ as samples from $\pi$ become available from the MCMC sampler. Finally, we point out at this stage that, although we restrict here our discussion to the IMH algorithm, the KL criterion can equally be used for other updates, such as the SRWM algorithm. The algorithm proposed by Haario, Saksman and Tamminen [19] is in this case a particular instance of the algorithm hereafter.

In order to allow for flexibility and the description of a general class of algorithms, we consider here mixtures of distributions in the exponential family for the auxiliary proposal distribution. More precisely, let $\Xi \subset \mathbb{R}^{n_\xi}$ and $\mathcal{Z} \subset \mathbb{R}^{n_z}$, for some integers $n_\xi$ and $n_z$, and define the following family of exponential probability densities (defined with respect to the product measure $\lambda^{\text{Leb}} \otimes \mu$ for some measure $\mu$ on $\mathcal{Z}$)

$$\mathcal{E}_c = \{f : f_\xi(x,z) = \exp\{-\psi(\xi) + \langle T(x,z), \phi(\xi) \rangle\}; \xi, x, z \in \Xi \times \mathsf{X} \times \mathcal{Z}\},$$

where $\psi : \Xi \to \mathbb{R}, \phi : \Xi \to \mathbb{R}^{n_\theta}$ and $T : \mathsf{X} \times \mathcal{Z} \to \mathbb{R}^{n_\theta}$. Let $\mathcal{E}$ denote the set of densities $\tilde{q}_\xi$ that are marginals of densities from $\mathcal{E}_c$, that is, such that for any $\xi, x \in \Xi \times \mathsf{X}$ we have

$$\tilde{q}_\xi(x) = \int_{\mathcal{Z}} f_\xi(x, z) \mu(dz). \tag{58}$$

This family of densities covers in particular finite mixtures of multivariate normal distributions. Here, the variable $z$ plays the role of the *label* of the class, which is not observed (see, e.g., [32]). Using standard missing data terminology, $f_\xi(x, z)$ is the *complete data likelihood* and $\tilde{q}_\xi$ is the associated *incomplete data likelihood*, which is the marginal of the complete data likelihood with respect to the class labels. When the number of observations is fixed, a classical approach to estimating the parameters of a mixture distribution consists of using the EM algorithm.

7.2. *Classical EM algorithm.* The classical EM algorithm is an iterative procedure which consists of two steps. Given $n$ independent samples $(X_1, \ldots, X_n)$ distributed marginally according to $\pi$: (1) *Expectation step*: calculate the conditional expectation of the complete data log-likelihood, given the observations and $\xi_k$ (the estimate of $\xi$ at iteration $k$)

$$\xi \mapsto Q(\xi, \xi_k) \stackrel{\text{def}}{=} \sum_{i=1}^n \mathbb{E}[\log(f_\xi(X_i, Z_i))|X_i, \xi_k].$$

(2) *Maximization step*: maximize the function $\xi \mapsto Q(\xi, \xi_k)$ with respect to $\xi$. The new estimate for $\xi$ is $\xi_{k+1} = \arg\max_{\xi \in \Xi} Q(\xi, \xi_k)$ (provided that it exists and is unique). The key property at the core of the EM algorithm is that the incomplete data likelihood $\prod_{i=1}^n \tilde{q}_{\xi_{k+1}}(X_i) \geq \prod_{i=1}^n \tilde{q}_{\xi_k}(X_i)$ is increased at each iteration, with equality if and only if $\xi_k$ is a stationary point (i.e., a local or global minimum or a saddle point). Under mild additional conditions (see, e.g., [33]), the EM algorithm therefore converges to stationary points of the marginal likelihood. Note that, when $n \to \infty$, under appropriate conditions, the renormalized incomplete data log-likelihood $n^{-1} \sum_{i=1}^n \log \tilde{q}_\xi(X_i)$ converges to $\mathbb{E}_\pi[\log \tilde{q}_\xi(X)]$, which is equal, up to a constant and a sign, to the Kullback–Leibler divergence between $\pi$ and $\tilde{q}_\xi$. In our particular setting, the classical batch form of the algorithm is as follows: first define for $\xi \in \Xi$ the conditional distribution

$$\nu_\xi(x,z) \stackrel{\text{def}}{=} \frac{f_\xi(x,z)}{\tilde{q}_\xi(x)}, \tag{59}$$

where $\tilde{q}_\xi$ is given by (58). Now, assuming that $\int_{\mathcal{Z}} |T(x,z)| \nu_\xi(x,z) \mu(dz) < \infty$, one can define for $x \in \mathsf{X}$ and $\xi \in \Xi$

$$\nu_\xi T(x) \stackrel{\text{def}}{=} \int_{\mathcal{Z}} T(x,z) \nu_\xi(x,z) \mu(dz), \tag{60}$$

and check that for $f_\xi \in \mathcal{E}_c$ and any $(\xi, \xi') \in \Xi \times \Xi$ that

$$\mathbb{E}\{\log(f_\xi(X_i, Z_i)) | X_i, \xi'\} = L(\nu_{\xi'} T(X_i); \xi),$$

where $L : \Theta \times \Xi \to \mathbb{R}$ is defined as

$$L(\theta; \xi) \stackrel{\text{def}}{=} -\psi(\xi) + \langle \theta, \phi(\xi) \rangle \qquad \text{with } \Theta \stackrel{\text{def}}{=} T(\mathsf{X}, \mathcal{Z}). \tag{61}$$

From this, one easily deduces that, for $n$ samples,

$$Q(\xi, \xi_k) = nL\left(\frac{1}{n} \sum_{i=1}^n \nu_{\xi_k} T(X_i); \xi\right).$$

Assuming now for simplicity that, for all $\theta \in \Theta$, the function $\xi \to L(\theta; \xi)$ reaches its maximum at a single point denoted by $\hat{\xi}(\theta)$ [i.e., $L(\theta; \hat{\xi}(\theta)) \geq L(\theta; \xi)$ for all $\xi \in \Xi$], the EM recursion can then be simply written as

$$\xi_{k+1} = \hat{\xi}\left(\frac{1}{n} \sum_{i=1}^n \nu_{\xi_k} T(X_i)\right).$$

The condition on the existence and uniqueness of $\hat{\xi}(\theta)$ is not restrictive. It is, for example, satisfied for finite mixtures of normal distributions. More sophisticated generalizations of the EM algorithm have been developed in order to deal with situations where this condition is not satisfied; see, for example, [25].

Our scenario differs from the classical setup above in two respects. First, the number of samples considered evolves with time, which requires that we estimate $\xi$ on the fly. Second, the samples $\{X_i\}$ are generated by a transition probability with invariant distribution $\pi$ and are therefore not independent. We address the first problem in Section 7.3 and the two problems simultaneously in Section 7.4 where we describe our particular adaptive MCMC algorithm.

7.3. *Sequential EM algorithm.* Sequential implementations of the EM algorithm for estimating the parameters of a mixture when the data are observed sequentially in time have been considered by several authors (see [32], Chapter 6, [5] and the references therein). The version presented here is, in many respects, a standard adaptation of these algorithms and consists of recursively and jointly estimating and maximizing with respect to $\xi$ the function

$$\theta(\xi) = \mathbb{E}_\pi[\log \tilde{q}_\xi(X)] = \pi\{\nu T_\xi(X)\},$$

which, as pointed out earlier, is the Kullback–Leibler divergence between $\pi$ and $\tilde{q}_\xi$, up to an additive constant and a sign. At iteration $k+1$, given an estimate $\theta_k$ of $\theta$ and $\xi_k = \hat{\xi}(\theta_k)$, sample $X_{k+1} \sim \pi$ and calculate

$$\begin{aligned}(62)\quad \theta_{k+1} &= (1-\gamma_{k+1})\theta_k + \gamma_{k+1}\nu_{\xi_k}T(X_{k+1}) \\ &= \theta_k + \gamma_{k+1}(\nu_{\xi_k}T(X_{k+1}) - \theta_k),\end{aligned}$$

where $\{\gamma_k\}$ is a sequence of stepsizes and $\gamma_k \in [0,1]$. This can be interpreted as a stochastic approximation algorithm $\theta_{k+1} = \theta_k + \gamma_{k+1}H(\theta_k, X_{k+1})$ with, for $\theta \in \Theta$,

$$(63)\qquad H(\theta, x) = \nu_{\hat{\xi}(\theta)}T(x) - \theta \quad \text{and} \quad h(\theta) = \pi(\nu_{\hat{\xi}(\theta)}T) - \theta.$$

At this stage, it is possible to introduce a set of simple conditions on the distributions in $\mathcal{E}_c$ that ensures the convergence of the sequence $\{\theta_k\}$ defined above. By convergence, we mean here that $\{\theta_k\}$ converges to the set of stationary points of the Kullback–Leibler divergence between $\pi$ and $\tilde{q}_{\hat{\xi}(\theta)}$, that is,

$$\mathcal{L} \stackrel{\text{def}}{=} \{\theta \in \Theta : \nabla w(\theta) = 0\},$$

where, for $\theta \in \Theta$

$$(64)\qquad w(\theta) = K(\pi \| \tilde{q}_{\hat{\xi}(\theta)}),$$

and $K$ and $\tilde{q}_\xi$ are defined in (57) and (58), respectively. It is worth noticing that these very same conditions will be used to prove the convergence of our adaptive MCMC algorithm:

(E1)   (i) The sets $\Xi$ and $\Theta$ are open subsets of $\mathbb{R}^{n_\xi}$ and $\mathbb{R}^{n_\theta}$, respectively. $\mathcal{Z}$ is a compact subset of $\mathbb{R}^{n_z}$.

(ii) For any $x \in \mathsf{X}$, $T(x) \stackrel{\text{def}}{=} \inf\{M : \mu(\{z : |T(x,z)| \geq M\}) = 0\} < \infty$.

(iii) The functions $\psi : \Xi \to \mathbb{R}$ and $\phi : \Xi \to \mathbb{R}^{n_\theta}$ are twice continuously differentiable on $\Xi$.

(iv) There exists a continuously differentiable function $\hat{\xi} : \Theta \to \Xi$ such that, for all $\theta, \xi \in \Theta \times \Xi$,

$$L(\theta; \hat{\xi}(\theta)) \geq L(\theta; \xi).$$

REMARK 5. For many models, the function $\xi \to L(\theta; \xi)$ admits a unique global maximum for any $\theta \in \Theta$, and the existence and differentiability of $\theta \to \hat{\xi}(\theta)$ follows from the implicit function theorem under mild regularity conditions.

(E2)  (i) The level sets $\{\theta \in \Theta, w(\theta) \leq M\}$ for $M > 0$ are compact.

(ii) The set $\mathcal{L} \stackrel{\text{def}}{=} \{\theta \in \Theta, \nabla w(\theta) = 0\}$ of stationary points is included in a compact subset of $\Theta$.

(iii) The closure of $w(\mathcal{L})$ has an empty interior.

REMARK 6. Assumption (E2) depends on both the properties of $\pi$ and $q_{\hat{\xi}(\theta)}$ and should therefore be checked on a case-by-case basis. Note, however, that (a) these assumptions are satisfied for finite mixtures of distributions in the exponential family under classical technical conditions on the parametrization beyond the scope of the present paper (see, among others [32], Chapter 6, and [5] for details) (b) the third assumption in (E2) can very often be checked using Sard's theorem.

We first prove here an intermediate proposition concerned with estimates of the variation $\tilde{q}_\xi - \tilde{q}_{\xi'}$ under (E2) in various senses. Note that most of these results are not used in this section, but will be useful in the following.

PROPOSITION 16. *Let $\{\tilde{q}_\xi, \xi \in \Xi\} \subset \mathcal{E}$ be a family of distributions satisfying* (E1). *Then, for any convex compact set $\mathcal{K} \subset \Xi$:*

1. *There exists a constant $C < \infty$ such that*

(65) $$\sup_{\xi \in \mathcal{K}} |\nabla_\xi \log \tilde{q}_\xi(x)| \leq C(1 + T(x)).$$

2. *For any $\xi, \xi', x \in \mathcal{K}^2 \times \mathsf{X}$ there exists a constant $C < \infty$ such that*

(66) $$|\tilde{q}_\xi(x) - \tilde{q}_{\xi'}(x)| < C|\xi - \xi'|(1 + T(x)) \sup_{\xi \in \mathcal{K}} \tilde{q}_\xi(x).$$

3. *For $W \to [1, \infty)$ such that $\sup_{\xi \in \mathcal{K}} \int_\mathsf{X} \tilde{q}_\xi(x)[1 + T(x)]W(x)\lambda^{\text{Leb}}(dx) < \infty$ and any $\xi, \xi' \in \mathcal{K}$, there exists a constant $C < \infty$ such that*

(67) $$\int_\mathsf{X} |\tilde{q}_\xi(x) - \tilde{q}_{\xi'}(x)|W(x)\lambda^{\text{Leb}}(dx) \leq C|\xi - \xi'|.$$

The proof appears in Appendix E. The key to establishing the convergence of the stochastic approximation procedure here consists of proving that $w(\theta) = K(\pi \| q_{\hat{\xi}(\theta)})$ plays the role of a Lyapunov function. This is hardly surprising as the algorithm aims at minimizing sequentially in time the incomplete data likelihood. More precisely, we have:

PROPOSITION 17. *Assume* (E1). *Then, for all* $\theta \in \Theta$, $\langle \nabla w(\theta), h(\theta) \rangle \leq 0$, *and*

$$\mathcal{L} = \{\theta \in \Theta : \langle \nabla w(\theta), h(\theta) \rangle = 0\} = \{\theta \in \Theta : \nabla w(\theta) = 0\}, \tag{68}$$

$$\hat{\xi}(\mathcal{L}) = \{\xi \in \Xi : \nabla_\xi K(\pi \| q_\xi) = 0\}, \tag{69}$$

*where* $\theta \mapsto h(\theta)$ *is given in* (63).

The proof appears in Appendix E. Another important result needed to prove convergence is the regularity of the field $\theta \mapsto H_\theta$. We have:

PROPOSITION 18. *Assume* (E1). *Then* $\{H_\theta, \theta \in \Theta\}$ *is* $(1+T)^2$-*Lipschitz, where* $H_\theta$ *is defined in* (63).

The proof appears in Appendix E. With this, and standard results on the convergence of SA, one may show that the SA procedure converges pointwise under (E1) and (E2).

7.4. *On-line EM for IMH adaptation.* We now consider the combination of the sequential EM algorithm described earlier with the IMH sampler. As we shall see in Proposition 20, using $\tilde{q}_{\hat{\xi}(\theta)}$ as a proposal distribution for the IMH transition is not sufficient to ensure the convergence of the algorithm, and it will be necessary to use a mixture of a *fixed* distribution $\zeta$ (which will not be updated during the successive iterations) and an *adaptive* component, here $\tilde{q}_{\hat{\xi}(\theta)}$. More precisely, we define the following family of parametrized IMH transition probabilities $\{P_\theta, \theta \in \Theta\}$: for $e \in (0,1]$ let $\zeta \in \mathcal{Q}_{e,\pi}$ (assumed nonempty) be a density which does not depend on $\theta \in \Theta$, let $\iota \in (0,1)$ and define the family of IMH transition probabilities

$$\mathcal{P}_{e,\zeta} \stackrel{\text{def}}{=} \{P_\theta \stackrel{\text{def}}{=} P_{q_\theta}^{\text{IMH}}, \theta \in \Theta\} \quad \text{with } \{q_\theta \stackrel{\text{def}}{=} (1-\iota)\tilde{q}_{\hat{\xi}(\theta)} + \iota\zeta, \theta \in \Theta\}. \tag{70}$$

The following properties on $\zeta$ and $\mathcal{E}_c$ will be required in order to ensure that $\mathcal{P}_{e,\zeta}$ satisfies (E1) and (E2):

(E3) (i) There exist $e > 0$ and $\zeta \in \mathcal{Q}_{e,\pi}$ such that, for any compact $\mathcal{K} \subset \Xi$,

$$\sup_{\xi \in \mathcal{K}} \inf \left\{ M : \lambda^{\text{Leb}}\left(\frac{\tilde{q}_\xi(x)(1+T(x))}{\zeta(x)} \geq M\right) = 0 \right\} < \infty. \tag{71}$$

(ii) There exists $W \to [1,\infty)$ such that, for any compact subset $\mathcal{K} \subset \Xi$,

$$\int_\mathsf{X} W(x)(1+T(x))\zeta(x)\lambda^{\text{Leb}}(dx)$$
$$+ \sup_{\xi \in \mathcal{K}} \int_\mathsf{X} W(x)(1+T(x))\tilde{q}_\xi(x)\lambda^{\text{Leb}}(dx) < \infty,$$

and $\sup_{x \in \mathsf{K}} W(x) < \infty$, where $\mathsf{K}$ is defined in Section 2.

REMARK 7. It is worth pointing out that the above choice for $q_\theta$ and the condition $\zeta \in \mathcal{Q}_{e,\pi}$ automatically ensure that $\{q_\theta, \theta \in \Theta\} \subset \mathcal{Q}_{\varepsilon,\pi}$ for $\varepsilon = e\iota$.

The *basic* version (see Section 2) of our algorithm now proceeds as follows: set $\theta_0 \in \Theta$, $\xi_0 = \hat{\xi}(\theta_0)$ and draw $X_0$ according to some initial distribution. At iteration $k+1$ for $k \geq 0$, draw $X_{k+1} \sim P_{\theta_k}(X_k, \cdot)$ where $P_\theta$ is given in (70). Compute $\theta_{k+1} = \theta_k + \gamma_{k+1}(\nu_{\xi_k} T(X_{k+1}) - \theta_k)$ and $\xi_{k+1} = \hat{\xi}(\theta_{k+1})$. We will study here the corresponding algorithm with reprojections which results in the homogeneous Markov chain $\{Z_k, k \geq 0\}$ as described in Section 2.

We now establish intermediate results about $\mathcal{P}_{e,\zeta}$ and $\{H_\theta, \theta \in \Theta\}$ which will lead to the proof that (A1)–(A3) are satisfied. We start with a general proposition about the properties of IMH transition probabilities, relevant to checking (A1).

PROPOSITION 19. *Let $V : \mathsf{X} \to [1, +\infty)$, and let $q \in \mathcal{Q}_{\varepsilon,\pi}$ for some $\varepsilon > 0$. Then:*

1. $\mathsf{X}$ *is a $(1, \varepsilon)$-small set with minorization distribution $\varphi = q$, and*

$$P_q^{\mathrm{IMH}} V(x) \leq (1-\varepsilon)V(x) + q(V), \qquad \text{where } q(V) = \int_{\mathsf{X}} q(x)V(x)\lambda^{\mathrm{Leb}}(dx).$$

2. *For any $f \in \mathcal{L}_V$, and any proposal distributions $q, q' \in \mathcal{Q}_{\varepsilon,\pi}$,*

$$(2\|f\|_V)^{-1} \|P_q^{\mathrm{IMH}} f - P_{q'}^{\mathrm{IMH}} f\|_V$$
$$(72) \qquad \leq \int_{\mathsf{X}} |q(x) - q'(x)| V(x) \lambda^{\mathrm{Leb}}(dx)$$
$$\qquad + [q(V) \vee q'(V)]((1 \wedge |1 - q^{-1} q'|_1) \vee (1 \wedge |1 - (q')^{-1} q|_1)).$$

The proof appears in Appendix F. In contrast with the SRWM, the $V$-norm $\|P_q^{\mathrm{IMH}} f - P_{q'}^{\mathrm{IMH}} f\|_V$ can be large, even in situations where $\int_{\mathsf{X}} |q(x) - q'(x)| \times V(x) \lambda^{\mathrm{Leb}}(dx)$ is small. This stems from the fact that the ratio of densities $q/q'$ enters the upper bound above. As we shall see in Proposition 20 below, this is what motivates our definition of the proposal distributions in (70) as a mixture of $\tilde{q}_\xi \in \mathcal{E}$ and a nonadaptive distribution $\zeta$ which satisfies (E3).

PROPOSITION 20. *Assume that the family of distributions $\{\tilde{q}_\xi, \xi \in \Xi\} \subset \mathcal{E}$ satisfies (E1) and (E3). Then, the family of transition kernels $\mathcal{P}_{e,\zeta}$ given in (70) satisfies (A1) with $\epsilon = e\iota$, $V = W$, $\varphi = \zeta$ and, if $W$ is bounded, then $\mathsf{C} = \mathsf{X}$, $\lambda = 0$, otherwise choose $\varepsilon \in (0, e\iota)$ such that $\mathsf{C} = \{x : V(x) < \varepsilon^{-1} \sup_{\theta \in \mathcal{K}} q_\theta(V)\}$ is such that $\zeta(\mathsf{C}) > 0$, and set $\lambda = 1 - e\iota + \varepsilon$.*

The proof appears in Appendix F. We are now in a position to present our final result:

THEOREM 21. *Let $\pi \in \mathcal{M}(\mathsf{X})$ and $\{\tilde{q}_\xi, \xi \in \Xi\} \subset \mathcal{E}$ be a family of distributions. Define $\Theta := T(\mathsf{X}, \mathcal{Z})$. Consider the following families of transition probabilities and functions:*

  (i) *$\{P_\theta, \theta \in \Theta\}$, as in (70), where $\zeta \in \mathcal{Q}_{e,\pi}$ for some $e > 0$, $\{\tilde{q}_\xi, \xi \in \Xi\}$ is further assumed to satisfy* (E1), (E3) *[with $V$ such that $T \in \mathcal{L}_{V^{\beta/2}}$ for some $\beta \in [0, 1)]$ and* (E2);
  (ii) *$\{H_\theta, \theta \in \Theta\}$ as in (63).*

*Let $\{\mathcal{K}_q, q \geq 0\}$ be a compact coverage of $\Theta$, let $\mathsf{K}$ be a compact set and let $\gamma = \{\gamma_k\}$ satisfy* (A5). *Consider the time-homogeneous Markov chain $\{Z_k\}$ on $\mathsf{Z}$ with transition probability $R$, as defined in Section 2. Then, for any $(x, \theta) \in \mathsf{K} \times \mathcal{K}_0$ and any $\alpha < 1 - \beta$:*

1. *For any $f \in \mathcal{L}_{V^\alpha}$,*
$$n^{-1} \sum_{k=1}^{n} (f(X_k) - \pi(f)) \to 0, \qquad \overline{\mathbb{P}}_\star\text{-}a.s.$$

2. *There $\overline{\mathbb{P}}_\star$-a.s. exists a random variable $\theta_\infty \in \{\theta \in \Theta : \nabla_\theta K(\pi \| \tilde{q}_{\hat{\xi}(\theta)}) = 0\}$ such that, provided that $\sigma(\theta_\infty, f) > 0$ and for any $f \in \mathcal{L}_{V^{\alpha/2}}$,*
$$\frac{1}{\sqrt{n}\sigma(\theta_\infty, f)} \sum_{k=1}^{n} (f(X_k) - \pi(f)) \xrightarrow{\mathcal{D}}_{\overline{\mathbb{P}}_\star} \mathcal{N}(0, 1),$$

*where $\sigma(\theta, f)$ is given as in (38).*

PROOF. The application of Propositions 17, 18 and 20 shows that (A1)–(A3) are satisfied, which, together with (E2) and (A5), implies Theorem 11. Then, we conclude by invoking Theorems 8 and 9. □

REMARK 8. It is worth noting that, provided $\pi \in \mathcal{M}(\mathsf{X})$ satisfies (M), the results of Propositions 12, 16, 17 and 18 proved in this paper easily allow one to establish a result similar to Theorem 21 for a generalization of the N-SRWM of [19] (described here in Section 1 and studied in Section 5) to the case where the proposal distribution belongs to $\mathcal{E}$, that is, when the proposal is a mixture of distributions.

## APPENDIX A: STABILITY OF THE INHOMOGENEOUS CHAIN

PROOF OF LEMMA 5. Under (A1) and (A2), we have, for $x, \theta \in \mathsf{X} \times \Theta$ and $k \geq 1$,
$$\mathbb{E}_{x,\theta}^{\boldsymbol{\rho}}[V(X_k)\mathbb{1}\{\sigma(\mathcal{K}) \geq k\}] = \mathbb{E}_{x,\theta}^{\boldsymbol{\rho}}[\mathbb{E}^{\boldsymbol{\rho}}[V(X_k)|\mathcal{F}_{k-1}]\mathbb{1}\{\sigma(\mathcal{K}) \geq k\}]$$

$$\leq \lambda \mathbb{E}^{\boldsymbol{\rho}}_{x,\theta}[V(X_{k-1})\mathbb{1}\{\sigma(\mathcal{K}) \geq k-1\}] + b.$$

Now, a straightforward induction leads to

$$\mathbb{E}^{\boldsymbol{\rho}}_{x,\theta}[V(X_k)\mathbb{1}\{\sigma(\mathcal{K}) \geq k\}] \leq \lambda^k V(x) + \frac{b}{1-\lambda},$$

which shows that there exists a constant $C$ [depending only on $\mathcal{K}$ and the constants appearing in (A1)] such that, for all $k \geq 0$,

(73) $$\mathbb{E}^{\boldsymbol{\rho}}_{x,\theta}[V(X_k)\mathbb{1}\{\sigma(\mathcal{K}) \geq k\}] \leq CV(x).$$

Now, for any $r \in [0,1]$, by Jensen's inequality, we have, for any $k \geq 0$,

(74) $\mathbb{E}^{\boldsymbol{\rho}}_{x,\theta}[V^r(X_k)\mathbb{1}\{\sigma(\mathcal{K}) \geq k\}] \leq (\mathbb{E}^{\boldsymbol{\rho}}_{x,\theta}[V(X_k)\mathbb{1}\{\sigma(\mathcal{K}) \geq k\}])^r \leq C^r V^r(x),$

showing (21). Similarly, again using Jensen's inequality,

$$\mathbb{E}^{\boldsymbol{\rho}}_{x,\theta}\left[\max_{1 \leq m \leq k}(a_m V(X_m))^r \mathbb{1}\{\sigma(\mathcal{K}) \geq m\}\right]$$
$$\leq \left(\mathbb{E}^{\boldsymbol{\rho}}_{x,\theta}\left[\max_{1 \leq m \leq k} a_m V(X_m)\mathbb{1}\{\sigma(\mathcal{K}) \geq m\}\right]\right)^r$$
$$\leq \left(\sum_{m=1}^{k} a_m \mathbb{E}^{\boldsymbol{\rho}}_{x,\theta}[V(X_m)\mathbb{1}\{\sigma(\mathcal{K}) \geq m\}]\right)^r \leq C^r \left(\sum_{m=1}^{k} a_m\right)^r V^r(x),$$

showing (22). Finally, since

$$\mathbb{1}\{\sigma(\mathcal{K}) \geq m\} \sum_{k=1}^{m} V^r(X_k) \leq \sum_{k=1}^{m} V^r(X_k)\mathbb{1}\{\sigma(\mathcal{K}) \geq k\},$$

we have

$$\mathbb{E}^{\boldsymbol{\rho}}_{x,\theta}\left[\max_{1 \leq m \leq n} \mathbb{1}\{\sigma(\mathcal{K}) \geq m\} \sum_{k=1}^{m} V^r(X_k)\right]$$
$$\leq \mathbb{E}^{\boldsymbol{\rho}}_{x,\theta}\left[\sum_{k=1}^{n} V^r(X_k)\mathbb{1}\{\sigma(\mathcal{K}) \geq k\}\right] \leq CnV^r(x),$$

showing (23). $\square$

The following proposition is a direct adaptation of Birnbaum and Marshall [10] inequality:

PROPOSITION 22. *Let $\{S_k, \mathcal{F}_k, k \geq 0\}$ be a submartingale, that is, $E(S_k|\mathcal{F}_{k-1}) \geq S_{k-1}$ a.e. Let $\{a_k > 0, 1 \leq k \leq n\}$ be a nonincreasing real-valued sequence. If $p \geq 1$ is such that $E|S_k|^p < \infty$ for $k \in \{1,\ldots,n\}$, then for*

$m \leq n$,

$$P\left\{\max_{m \leq k \leq n} a_k |S_k| \geq 1\right\} \leq \sum_{k=m}^{n-1} (a_k^p - a_{k+1}^p) E|S_k|^p + a_n^p E|S_n|^p.$$

## APPENDIX B: PROOF OF PROPOSITIONS 3 AND 4

In the sequel, $C$ is a generic constant, which may take different values upon each appearance.

PROOF OF PROPOSITION 3. Let $\mathcal{K} \subset \Theta$ be a compact set and let $r \in [0,1]$. For any $(\theta, \theta') \in \mathcal{K} \times \mathcal{K}$ and $f \in \mathcal{L}(V^r)$,

$$P_\theta^n f - P_{\theta'}^n f = \sum_{j=0}^{n-1} P_\theta^j (P_\theta - P_{\theta'}) P_{\theta'}^{n-j-1} f$$

$$= \sum_{j=0}^{n-1} (P_\theta^j - \pi)(P_\theta - P_{\theta'})(P_{\theta'}^{n-j-1} f - \pi(f)),$$

where we have used the fact that $\pi P_\theta = \pi P_{\theta'} = \pi$ for any $\theta, \theta'$. Theorem 2 shows that there exists a constant $C$ and $\rho \in (0,1)$ such that, for any $\theta \in \mathcal{K}$, $l \geq 0$ and any $f \in \mathcal{L}_{V^r}$,

(75) $$\|P_\theta^l f - \pi(f)\|_{V^r} \leq C \|f\|_{V^r} \rho^l.$$

Under assumption (A2), for any $(\theta, \theta') \in \mathcal{K} \times \mathcal{K}$, $l \geq 0$ and any $f \in \mathcal{L}_{V^r}$,

$$\|(P_\theta^j - \pi)(P_\theta - P_{\theta'})(P_{\theta'}^{n-j-1} f - \pi(f))\|_{V^r}$$

$$\leq C \rho^j \|(P_\theta - P_{\theta'})(P_{\theta'}^{n-j-1} f - \pi(f))\|_{V^r}$$

$$\leq C \rho^j |\theta - \theta'| \|P_{\theta'}^{n-j-1} f - \pi(f)\|_{V^r}$$

$$\leq C |\theta - \theta'| \|f\|_{V^r} \rho^n,$$

showing that there exists a constant $C < \infty$ such that, for any $(\theta, \theta') \in \mathcal{K} \times \mathcal{K}$ and $f \in \mathcal{L}_{V^r}$,

(76) $$\|P_\theta^n f - P_{\theta'}^n f\|_{V^r} \leq C n \rho^n |\theta - \theta'| \|f\|_{V^r}.$$

Now consider $\{f_\theta, \theta \in \Theta\}$, a family of $V^r$-Lipschitz functions. From (75), for any $\theta \in \mathcal{K}$, $\sum_{k=0}^\infty |P_\theta^k f_\theta - \pi(f_\theta)| < \infty$ and $\hat{f}_\theta \stackrel{\text{def}}{=} \sum_{k=0}^\infty (P_\theta^k f_\theta - \pi(f_\theta))$ belongs to $\mathcal{L}_{V^r}$. Now we consider the difference

$$\hat{f}_\theta - \hat{f}_{\theta'} = \sum_{k=0}^\infty (P_\theta^k f_\theta - \pi(f_\theta)) - \sum_{k=0}^\infty (P_{\theta'}^k f_{\theta'} - \pi(f_{\theta'}))$$

$$= \sum_{k=0}^\infty (P_\theta^k f_\theta - P_{\theta'}^k f_\theta) - \sum_{k=0}^\infty (P_{\theta'}^k (f_{\theta'} - f_\theta) - \pi(f_{\theta'} - f_\theta)),$$

which using (75) and (76) shows that

$$\|\hat{f}_\theta - \hat{f}_{\theta'}\|_{V^r} \leq C|\theta - \theta'|\left(\sum_{k=0}^{\infty} k\rho^k\right)\|f_\theta\|_{V^r} + C\left(\sum_{k=0}^{\infty} \rho^k\right)\|f_\theta - f_{\theta'}\|_{V^r},$$

and we conclude by using the fact that $\{f_\theta, \theta \in \Theta\}$ is a $V^r$-Lipschitz family of functions. Using the same arguments one can prove a similar bound for $\|P_\theta \hat{f}_\theta - P_{\theta'}\hat{f}_{\theta'}\|_{V^r}$. $\square$

PROOF OF PROPOSITION 4. For simplicity, we set $\sigma := \sigma(\mathcal{K})$ and, in what follows, $C$ is a finite constant whose value might change upon each appearance. Let $x, \theta \in \mathsf{X} \times \Theta$. For $k \geq k_0$, we introduce the following decomposition:

$$|\mathbb{E}_{x,\theta}^{\bar{\rho}}\{(f(X_k) - \pi(f))\mathbb{1}(\sigma \geq k)\}|$$

$$\leq |\mathbb{E}_{x,\theta}^{\bar{\rho}}\{(f(X_k) - P_{\theta_{k-n(k)}}^{n(k)} f(X_{k-n(k)}))\mathbb{1}(\sigma \geq k)\}|$$

$$+ |\mathbb{E}_{x,\theta}^{\bar{\rho}}\{(P_{\theta_{k-n(k)}}^{n(k)} f(X_{k-n(k)}) - \pi(f))\mathbb{1}(\sigma \geq k)\}|.$$

By Theorem 2, the last term is bounded by $C\|f\|_{V^{1-\beta}} V^{1-\beta}(x)\rho^{n(k)} \leq C\bar{\rho}^{-1}\|f\|_{V^{1-\beta}} \bar{\rho}_k V(x)$. We consider the first term and use the following new decomposition of this bias term (cf. [19]),

$$|\mathbb{E}_{x,\theta}^{\bar{\rho}}\{(f(X_k) - P_{\theta_{k-n(k)}}^{n(k)} f(X_{k-n(k)}))\mathbb{1}(\sigma \geq k)\}|$$

$$\leq \left|\sum_{j=2}^{n(k)} \mathbb{E}_{x,\theta}^{\bar{\rho}}\{(P_{\theta_{k-j+1}}^{j-1} f(X_{k-j+1}) - P_{\theta_{k-j}}^{j} f(X_{k-j}))\mathbb{1}(\sigma > k - j + 1)\}\right|$$

$$\leq \left|\sum_{j=2}^{n(k)} \mathbb{E}_{x,\theta}^{\bar{\rho}}\{\mathbb{E}_{x,\theta}^{\bar{\rho}}\{P_{\theta_{k-j+1}}^{j-1} f(X_{k-j+1}) - P_{\theta_{k-j}}^{j-1} f(X_{k-j+1})|\mathcal{F}_{k-j}\} \right.$$

$$\left. \times \mathbb{1}(\sigma > k - j + 1)\}\right|$$

$$\leq C\|f\|_{V^{1-\beta}} \sum_{j=2}^{n(k)} (j-1)\rho^{j-1} |\mathbb{E}_{x,\theta}^{\bar{\rho}}\{|\theta_{k-j+1} - \theta_{k-j}|V^{1-\beta}(X_{k-j+1})$$

$$\times \mathbb{1}(\sigma > k - j + 1)\}|$$

$$\leq C\|f\|_{V^{1-\beta}} \sum_{j=1}^{n(k)-1} j\rho^j \bar{\rho}_{k-j} \mathbb{E}_{x,\theta}^{\bar{\rho}}\{V(X_{k-j})\mathbb{1}(\sigma \geq k - j)\}$$

$$\leq C\|f\|_{V^{1-\beta}} V(x)\bar{\rho}_{k-n(k)+1},$$

where we have successively used (76) with $r = 1 - \beta$, (A3), the fact that $\bar{\rho}$ is assumed nonincreasing and (74). We conclude with the additional condition on $\bar{\rho}$. □

## APPENDIX C: PROOF OF PROPOSITION 12

For any $x \in \mathsf{X}$, define the acceptance region $\mathsf{A}(x) = \{z \in \mathsf{X}; \pi(x+z) \geq \pi(x)\}$ and the rejection region $\mathsf{R}(x) = \{z \in \mathsf{X}; \pi(x+z) < \pi(x)\}$. From the definition (47) of $\mathcal{Q}_{a,b}$ ([31], Theorem 2.2) applies for any $q \in \mathcal{Q}_{a,b}$ and we can conclude that (48) is satisfied. Noting that the two sets $\mathsf{A}(x)$ and $\mathsf{R}(x)$ do not depend on the proposal distribution $q$ and using the conclusion of the proof of Theorem 4.3 of [22], we have

$$\inf_{q \in \mathcal{Q}_{a,b}} \liminf_{|x| \to +\infty} \int_{\mathsf{A}(x)} q(z) \lambda^{\mathrm{Leb}}(dz) > 0,$$

so that, from the conclusion of the proof of Theorem 4.1 of [22],

$$\sup_{q \in \mathcal{Q}_{a,b}} \limsup_{|x| \to +\infty} \frac{P_q^{\mathrm{SRW}} V(x)}{V(x)} = 1 - \inf_{q \in \mathcal{Q}_{a,b}} \liminf_{|x| \to +\infty} \int_{\mathsf{A}(x)} q(z) \lambda^{\mathrm{Leb}}(dz) < 1,$$

which proves (49). Finally, for any $q \in \mathcal{K}_{a,b}$,

$$\frac{P_q^{\mathrm{SRW}} V(x)}{V(x)} = \int_{\mathsf{A}(x)} \frac{\pi(x+z)^{-\eta}}{\pi(x)^{-\eta}} q(z) \lambda^{\mathrm{Leb}}(dz)$$

$$+ \int_{\mathsf{R}(x)} \left(1 - \frac{\pi(x+z)}{\pi(x)} + \frac{\pi(x+z)^{1-\eta}}{\pi(x)^{1-\eta}}\right) q(z) \lambda^{\mathrm{Leb}}(dz)$$

$$\leq \sup_{0 \leq u \leq 1} (1 - u + u^{1-\eta}),$$

which proves (50). Now notice that

$$P_q^{\mathrm{SRW}} f(x) - P_{q'}^{\mathrm{SRW}} f(x) = \int_{\mathsf{X}} \alpha(x, x+z)(q(z) - q'(z))f(x+z)\lambda^{\mathrm{Leb}}(dz)$$

$$+ f(x) \int_{\mathsf{X}} \alpha(x, x+z)(q'(z) - q(z))\lambda^{\mathrm{Leb}}(dz).$$

We therefore focus, for $r \in [0,1]$ and $f \in \mathcal{L}_{V^r}$, on the term

$$\frac{|\int_{\mathsf{X}} \alpha(x, x+z)(q(z) - q'(z))f(x+z)\lambda^{\mathrm{Leb}}(dz)|}{\|f\|_{V^r} V^r(x)}$$

$$\leq \frac{\int_{\mathsf{X}} \alpha(x, x+z)|q(z) - q'(z)|V^r(x+z)\lambda^{\mathrm{Leb}}(dz)}{V^r(x)}$$

$$\leq \int_{\mathsf{A}(x)} \frac{\pi(x+z)^{-r\eta}}{\pi(x)^{-r\eta}} |q(z) - q'(z)|\lambda^{\mathrm{Leb}}(dz)$$

$$+ \int_{\mathsf{R}(x)} \frac{\pi(x+z)^{1-r\eta}}{\pi(x)^{1-r\eta}} |q(z) - q'(z)| \lambda^{\text{Leb}}(dz)$$

$$\leq \int_{\mathsf{X}} |q(z) - q'(z)| \lambda^{\text{Leb}}(dz).$$

We now conclude that, for any $x \in \mathsf{X}$ and any $f \in \mathcal{L}_{V^r}$,

$$\frac{|P_q^{\text{SRW}} f(x) - P_{q'}^{\text{SRW}} f(x)|}{V^r(x)} \leq 2\|f\|_{V^r} \int_{\mathsf{X}} |q(z) - q'(z)| \lambda^{\text{Leb}}(dz).$$

## APPENDIX D: PROOF OF LEMMA 13

For notational simplicity, we write $q_\Gamma$ for $\mathcal{N}(0, \Gamma)$. We have

$$\int_{\mathsf{X}} |q_\Gamma(z) - q_{\Gamma'}(z)| \lambda^{\text{Leb}}(dz) = \int_{\mathsf{X}} \left| \int_0^1 \frac{d}{dv} q_{\Gamma + v(\Gamma' - \Gamma)}(z) \lambda^{\text{Leb}}(dv) \right| \lambda^{\text{Leb}}(dz),$$

and let $\Gamma_v = \Gamma + v(\Gamma' - \Gamma)$, so that

$$\frac{d}{dv} \log q_{\Gamma + v(\Gamma' - \Gamma)}(z) = -\frac{1}{2} \text{Tr}[\Gamma_v^{-1}(\Gamma' - \Gamma) + \Gamma_v^{-1} z z^{\mathsf{T}} \Gamma_v^{-1}(\Gamma' - \Gamma)],$$

and consequently

$$\int_{\mathsf{X}} \left| \int_0^1 \frac{d}{dv} q_{\Gamma + v(\Gamma' - \Gamma)}(z) \lambda^{\text{Leb}}(dv) \right| \lambda^{\text{Leb}}(dz) \leq |\Gamma' - \Gamma| \int_0^1 |\Gamma_v^{-1}| \lambda^{\text{Leb}}(dv)$$

$$\leq \frac{n_x}{\lambda_{\min}(\mathcal{K})} |\Gamma' - \Gamma|,$$

where we have used the following inequality:

$$|\text{Tr}[\Gamma_v^{-1} z z^{\mathsf{T}} \Gamma_v^{-1}(\Gamma' - \Gamma)]| \leq |\Gamma' - \Gamma| \text{Tr}[\Gamma_v^{-1} \Gamma_v^{-1} z z^{\mathsf{T}}].$$

## APPENDIX E: PROOFS OF PROPOSITIONS 16, 17 AND 18

Hereafter, for a scalar function $s$, $\nabla s$ is a column vector, and for a (column) vector-valued function $v$ with scalar entries $v_1, v_2, \ldots$, we use the convention that $\nabla v$ is the matrix with $\nabla v_i$ as its $i$th column.

PROOF OF PROPOSITION 16. We first note that, from *Fisher's identity*, we have

$$\forall \xi \in \Xi \quad \nabla_\xi \log \tilde{q}_\xi(x) = \int_{\mathcal{Z}} \nabla_\xi \log f_\xi(x, z) \nu_\xi(x, z) \mu(dz)$$

$$= -\nabla_\xi \psi(\xi) + \nabla_\xi \phi(\xi) \nu_\xi T(x)$$

and, from (E1), we conclude that (65) holds. Equation (66) is a direct consequence of (65). Now we prove (67). With $\Delta_\psi \stackrel{\text{def}}{=} \psi(\xi') - \psi(\xi)$ and

$$\Delta_\phi \stackrel{\text{def}}{=} \phi(\xi') - \phi(\xi) \text{ for } \xi, \xi' \in \Xi,$$

$$\int_{\mathsf{X}} |\tilde{q}_\xi(x) - \tilde{q}_{\xi'}(x)| W(x) \lambda^{\text{Leb}}(dx)$$

$$= \int_{\mathsf{X}} \left| \int_0^1 \int_{\mathcal{Z}} [\Delta_\psi - \langle T(x,z), \Delta_\phi \rangle] f_\xi^v(x,z) f_{\xi'}^{1-v}(x,z) \mu(dz) \lambda^{\text{Leb}}(dv) \right|$$
$$\times W(x) \lambda^{\text{Leb}}(dx)$$

$$\leq |\Delta_\psi| \int_{\mathsf{X}} \int_0^1 \int_{\mathcal{Z}} f_\xi^v(x,z) f_{\xi'}^{1-v}(x,z) W(x) \mu(dz) \lambda^{\text{Leb}}(dv\, dx)$$
$$+ |\Delta_\phi| \int_{\mathsf{X}} \int_0^1 T(x) \int_{\mathcal{Z}} f_\xi^v(x,z) f_{\xi'}^{1-v}(x,z) W(x) \mu(dz) \lambda^{\text{Leb}}(dv\, dx)$$

$$\leq |\Delta_\psi| \int_0^1 \left[ \int_{\mathsf{X}} W(x) \tilde{q}_\xi(x) \lambda^{\text{Leb}}(dx) \right]^v$$
$$\times \left[ \int_{\mathsf{X}} W(x) \tilde{q}_{\xi'}(x) \lambda^{\text{Leb}}(dx) \right]^{1-v} \lambda^{\text{Leb}}(dv)$$
$$+ |\Delta_\phi| \int_0^1 \left[ \int_{\mathsf{X}} T(x) W(x) \tilde{q}_\xi(x) \lambda^{\text{Leb}}(dx) \right]^v$$
$$\times \left[ \int_{\mathsf{X}} T(x) W(x) \tilde{q}_{\xi'}(x) \lambda^{\text{Leb}}(dx) \right]^{1-v} \lambda^{\text{Leb}}(dv),$$

and we conclude by invoking the assumptions on $W$, $\phi$ and $\psi$. $\square$

PROOF OF PROPOSITION 17. We first note that, from *Fisher's identity*, we have

$$\forall \xi \in \Xi \quad \nabla_\xi \log \tilde{q}_\xi(x) = \int_{\mathcal{Z}} \nabla_\xi \log f_\xi(x,z) \nu_\xi(x,z) \mu(dz)$$
$$= -\nabla_\xi \psi(\xi) + \nabla_\xi \phi(\xi) \nu_\xi T(x).$$

From (65) and (E1), we may derive under the sum sign to show that

$$\nabla_\xi \int_{\mathsf{X}} \pi(x) \log \tilde{q}_\xi(x) \lambda^{\text{Leb}}(dx) = \int_{\mathsf{X}} \pi(x) \nabla_\xi \log \tilde{q}_\xi(x) \lambda^{\text{Leb}}(dx)$$
$$= -\nabla_\xi \psi(\xi) + \nabla_\xi \phi(\xi) \pi(\nu_\xi T),$$

and thus, by the chain rule of derivations,

$$\nabla_\theta w(\theta) = -\nabla_\theta \hat{\xi}(\theta)(-\nabla_\xi \psi(\hat{\xi}(\theta)) + \nabla_\xi \phi(\hat{\xi}(\theta)) \pi(\nu_{\hat{\xi}(\theta)} T)).$$

For any $\theta \in \Theta$, $\hat{\xi}(\theta)$ is a stationary point of the mapping $\xi \to L(\theta, \xi)$ and, thus,

$$\nabla_\xi L(\theta, \hat{\xi}(\theta)) = -\nabla_\xi \psi(\hat{\xi}(\theta)) + \nabla_\xi \phi(\hat{\xi}(\theta))\theta = 0.$$

Consequently, (63) implies that $\nabla_\theta w(\theta) = -\nabla_\theta \hat{\xi}(\theta) \nabla_\xi \phi(\hat{\xi}(\theta)) h(\theta)$. We also notice that $\nabla_\theta \nabla_\xi L(\theta, \xi) = \nabla_\xi \phi(\xi)^{\mathrm{T}}$. Differentiation with respect to $\theta$ of the mapping $\theta \mapsto \nabla_\xi L(\theta, \hat{\xi}(\theta))$ yields

$$\nabla_\theta \nabla_\xi L(\theta, \hat{\xi}(\theta)) = \nabla_\xi \phi(\hat{\xi}(\theta))^{\mathrm{T}} + \nabla_\theta \hat{\xi}(\theta) \nabla_\xi^2 L(\theta, \hat{\xi}(\theta)) = 0.$$

We finally have

$$\langle \nabla_\theta w(\theta), h(\theta) \rangle = h(\theta)^{\mathrm{T}} \nabla_\theta \hat{\xi}(\theta) \nabla_\xi^2 L(\theta, \hat{\xi}(\theta)) (\nabla_\theta \hat{\xi}(\theta))^{\mathrm{T}} h(\theta),$$

which concludes the proof, since, under (E1), $\nabla_\xi^2 L(\theta, \hat{\xi}(\theta)) \leq 0$, for any $\theta \in \Theta$. $\square$

PROOF OF PROPOSITION 18. For any $x \in \mathsf{X}$,

$$(77) \quad |H_\theta(x) - H_{\theta'}(x)| \leq T(x) \int_{\mathcal{Z}} |\nu_{\hat{\xi}(\theta)}(x,z) - \nu_{\hat{\xi}(\theta')}(x,z)| \mu(dz) + |\theta' - \theta|.$$

From Proposition 16 one has that, for any compact set $\mathcal{K} \subset \Xi$, there exists a constant $C$ such that, for all $\xi, z \in \mathcal{K} \times \mathcal{Z}$,

$$|\nabla_\xi \log f_\xi(x,z)| \leq C(1 + T(x)) \quad \text{and} \quad |\nabla_\xi \log q(x;\xi)| \leq C(1 + T(x)).$$

Thus,

$$|\nabla_\xi \log \nu_\xi(x,z)| \leq |\nabla_\xi \log f_\xi(x,z)| + |\nabla_\xi \log q_\xi(x)| \leq 2C(1 + T(x)).$$

Hence, for all $\xi, \xi' \in \mathcal{K}$ and $z \in \mathcal{Z}$,

$$|\nu_\xi(x,z) - \nu_{\xi'}(x,z)| \leq 2C(1 + T(x))|\xi - \xi'|,$$

which, together with equation (77), concludes the proof. $\square$

## APPENDIX F: PROOFS OF PROPOSITIONS 19 AND 20

PROOF OF PROPOSITION 19. The minorization condition is a classical result; see [26]. Now notice that

$$\begin{aligned}
P_q^{\mathrm{IMH}} V(x) &= \int_{\mathsf{X}} \alpha_q(x,y) V(y) q(y) \lambda^{\mathrm{Leb}}(dy) \\
&\quad + V(x) \int_{\mathsf{X}} [1 - \alpha_q(x,y)] q(y) \lambda^{\mathrm{Leb}}(dy) \\
&\leq \left(1 - \int_{\mathsf{X}} \left(\frac{q(x)}{\pi(x)} \wedge \frac{q(y)}{\pi(y)}\right) \pi(y) \lambda^{\mathrm{Leb}}(dy)\right) V(x) + q(V),
\end{aligned}$$

where $\alpha_q$ is given in (55). The drift condition follows.

42    C. ANDRIEU AND E. MOULINES

From the definition of the transition probability, and for any $f \in \mathcal{L}_V$,

$$|P_q^{\mathrm{IMH}} f(x) - P_{q'}^{\mathrm{IMH}} f(x)|$$

$$\leq \|f\|_V \left\{ \int_{\mathsf{X}} |\alpha_q(x,y)q(y) - \alpha_{q'}(x,y)q'(y)| V(y) \lambda^{\mathrm{Leb}}(dy) \right.$$

$$\left. + V(x) \int_{\mathsf{X}} |\alpha_{q'}(x,y)q'(y) - \alpha_q(x,y)q(y)| \lambda^{\mathrm{Leb}}(dy) \right\}$$

$$\leq 2\|f\|_V V(x) \int_{\mathsf{X}} |\alpha_q(x,y)q(y) - \alpha_{q'}(x,y)q'(y)| V(y) \lambda^{\mathrm{Leb}}(dy).$$

We therefore bound

$$I = \int_{\mathsf{X}} \left| \frac{q(y)}{\pi(y)} \wedge \frac{q(x)}{\pi(x)} - \frac{q'(y)}{\pi(y)} \wedge \frac{q'(x)}{\pi(x)} \right| \pi(y) V(y) \lambda^{\mathrm{Leb}}(dy).$$

We introduce the following sets:

$$\mathsf{A}_q(x) = \left\{ y : \frac{q(y)}{\pi(y)} \leq \frac{q(x)}{\pi(x)} \right\} \quad \text{and} \quad \mathsf{B}_q(x) = \left\{ y : \frac{q(y)}{\pi(y)} \leq \frac{q'(x)}{\pi(x)} \right\},$$

and note that the following inequalities hold:

(78)
$$\forall y \in \mathsf{A}_{q'}^c(x) \cap \mathsf{A}_q^c(x) \qquad \pi(y) < \frac{\pi(x)}{q(x)} q(y) \wedge \frac{\pi(x)}{q'(x)} q'(y) \quad \text{and}$$

$$\forall y \in \mathsf{A}_{q'}^c(x) \cap \mathsf{B}_q^c(x) \qquad \pi(y) < \frac{\pi(x)}{q'(x)} (q'(y) \wedge q(y)).$$

We now decompose $I$ into four terms $I \stackrel{\mathrm{def}}{=} \sum_{i=1}^{4} I_i$, where

$$I = \int_{\mathsf{A}_q \cap \mathsf{A}_{q'}} \left| \frac{q(y)}{\pi(y)} - \frac{q'(y)}{\pi(y)} \right| \pi(y) V(y) \lambda^{\mathrm{Leb}}(dy)$$

$$+ \int_{\mathsf{A}_q^c \cap \mathsf{A}_{q'}^c} \left| \frac{q(x)}{\pi(x)} - \frac{q'(x)}{\pi(x)} \right| \pi(y) V(y) \lambda^{\mathrm{Leb}}(dy)$$

$$+ \int_{\mathsf{A}_q \cap \mathsf{A}_{q'}^c} \left| \frac{q(y)}{\pi(y)} - \frac{q'(x)}{\pi(x)} \right| \pi(y) V(y) \lambda^{\mathrm{Leb}}(dy)$$

$$+ \int_{\mathsf{A}_q^c \cap \mathsf{A}_{q'}} \left| \frac{q(x)}{\pi(x)} - \frac{q'(y)}{\pi(y)} \right| \pi(y) V(y) \lambda^{\mathrm{Leb}}(dy).$$

Here we have dropped $x$ in the set notation for simplicity. We now determine bounds for $I_i$, $i = 2, 3$. Notice that, since $y \in \mathsf{A}_q^c \cap \mathsf{A}_{q'}^c$,

$$I_2 \leq \left\{ \left| 1 - \frac{q'(x)}{q(x)} \right| \int_{\mathsf{A}_q^c \cap \mathsf{A}_{q'}^c} V(y) q(y) \lambda^{\mathrm{Leb}}(dy) \right\}$$

$$\wedge \left\{ \left|1 - \frac{q(x)}{q'(x)}\right| \int_{\mathsf{A}_q^c \cap \mathsf{A}_{q'}^c} V(y) q'(y) \lambda^{\text{Leb}}(dy) \right\}$$

$$\leq \left\{ \left|1 - \frac{q'(x)}{q(x)}\right| \wedge \left|1 - \frac{q(x)}{q'(x)}\right| \right\}$$

$$\times \left\{ \int_{\mathsf{A}_q^c \cap \mathsf{A}_{q'}^c} V(y) q(y) \lambda^{\text{Leb}}(dy) \vee \int_{\mathsf{A}_q^c \cap \mathsf{A}_{q'}^c} V(y) q'(y) \lambda^{\text{Leb}}(dy) \right\},$$

and it can easily be checked that

$$\left|1 - \frac{q'}{q}\right| \wedge \left|1 - \frac{q}{q'}\right| \leq \left\{1 \wedge \left|1 - \frac{q'}{q}\right|\right\} \vee \left\{1 \wedge \left|1 - \frac{q}{q'}\right|\right\}.$$

The term $I_3$ can be bounded as follows:

$$I_3 \leq \left\{ \int_{\mathsf{A}_q \cap \mathsf{A}_{q'}^c \cap \mathsf{B}_q^c} q(y) V(y) \lambda^{\text{Leb}}(dy) \right\}$$

$$\wedge \left\{ \left( \frac{q(x)}{\pi(x)} - \frac{q'(x)}{\pi(x)} \right) \int_{\mathsf{A}_q \cap \mathsf{A}_{q'}^c \cap \mathsf{B}_q^c} V(y) \pi(y) \lambda^{\text{Leb}}(dy) \right\}$$

$$+ \int_{\mathsf{A}_q \cap \mathsf{A}_{q'}^c \cap \mathsf{B}_q} \left( \frac{q'(y)}{\pi(y)} - \frac{q(y)}{\pi(y)} \right) V(y) \pi(y) \lambda^{\text{Leb}}(dy),$$

and using (78) we find that

$$I_3 \leq \left\{ 1 \wedge \left( \frac{q(x)}{q'(x)} - 1 \right) \right\} \int_{\mathsf{A}_q \cap \mathsf{A}_{q'}^c \cap \mathsf{B}_q^c} q(y) V(y) \lambda^{\text{Leb}}(dy)$$

$$+ \int_{\mathsf{A}_q \cap \mathsf{A}_{q'}^c \cap \mathsf{B}_q} |q'(y) - q(y)| V(y) \lambda^{\text{Leb}}(dy).$$

The bound for $I_4$ follows from that of $I_3$ by swapping $q$ and $q'$. □

PROOF OF PROPOSITION 20. The first claim follows directly from Proposition 19 and the assumptions. Now denote

$$\Upsilon_{\xi,\xi',\alpha}(x) \stackrel{\text{def}}{=} \frac{(1-\alpha)\tilde{q}_\xi(x) + \alpha\zeta(x)}{(1-\alpha)\tilde{q}_{\xi'}(x) + \alpha\zeta(x)} = 1 + \frac{\tilde{q}_\xi(x) - \tilde{q}_{\xi'}(x)}{\zeta(x)[\tilde{q}_{\xi'}(x)/\zeta(x) + \alpha/(1-\alpha)]}.$$

Therefore, from (65), for any convex compact set $\mathcal{K} \subset \Xi$, there exists $C < \infty$ such that, for any $\xi, \xi', x \in \mathcal{K}^2 \times \mathsf{X}$,

$$|1 - \Upsilon_{\xi,\xi',\alpha}(x)| \leq \frac{1-\alpha}{\alpha} \frac{|\tilde{q}_\xi(x) - \tilde{q}_{\xi'}(x)|}{\zeta(x)} \leq C|\xi - \xi'| \frac{\sup_{\xi \in \mathcal{K}} \tilde{q}_\xi(x)(1 + T(x))}{\zeta(x)},$$

which, with (71), implies that, for all $\xi, \xi' \in \mathcal{K}$ and for $\lambda^{\text{Leb}}$-almost all $x$, there exists $C < \infty$ such that

$$(1 \wedge |1 - \Upsilon_{\xi,\xi',\alpha}(x)|) \vee (1 \wedge |1 - \Upsilon_{\xi',\xi,\alpha}(x)|) \leq C|\xi - \xi'|.$$

Now, as a direct consequence of equation (67), one can show that there exists $C$ such that for any $\xi, \xi' \in \mathcal{K}$ and $r \in [0,1]$,

$$\int_{\mathsf{X}} |\tilde{q}_\xi(x) - \tilde{q}_{\xi'}(x)| W(x)^r \lambda^{\text{Leb}}(dx) \leq C|\xi - \xi'|.$$

The proof is concluded by application of Proposition 19. □

**Acknowledgments.** The authors would like to thank Randal Douc, David Leslie, the anonymous reviewers and in particular Gersende Fort for very helpful comments that have helped improve an earlier version of the paper.

The authors would like to thank the Royal Society/CNRS partnership and the Royal Society International Exchange visit program, 2004–2006. Part of this work was presented at the First Cape Cod Monte Carlo Workshop, September 13–14, 2002.

## REFERENCES

[1] ANDRIEU, C. (2004). Discussion on the meeting on statistical approaches to inverse problems, December 2003. *J. R. Statist. Soc. Ser. B Stat. Methodol.* **66** 633–634.

[2] ANDRIEU, C. and MOULINES, E. (2003). On the ergodicity property of some adaptive MCMC algorithms. Technical report, Univ. Bristol.

[3] ANDRIEU, C., MOULINES, E. and PRIOURET, P. (2005). Stability of stochastic approximation under verifiable conditions. *SIAM J. Control Optim.* **44** 283–312. MR2177157

[4] ANDRIEU, C. and ROBERT, C. (2001). Controlled MCMC for optimal sampling. Cahiers du Cérémade 0125.

[5] ARCIDIACONO, P. and JONES, J. B. (2003). Finite mixture distributions, sequential likelihood and the EM algorithm. *Econometrica* **71** 933–946. MR1983232

[6] ATCHADÉ, Y. and ROSENTHAL, J. (2003). On adaptive Markov chain Monte Carlo algorithms. Technical report.

[7] BARTUSEK, J. (2000). Stochastic approximation and optimization of Markov chains. Ph.D. thesis.

[8] BAXENDALE, P. H. (2005). Renewal theory and computable convergence rates for geometrically ergodic Markov chains. *Ann. Appl. Probab.* **15** 700–738. MR2114987

[9] BENVENISTE, A., MÉTIVIER, M. and PRIOURET, P. (1990). *Adaptive Algorithms and Stochastic Approximations.* Springer, New York.

[10] BIRNBAUM, Z. W. and MARSHALL, A. W. (1961). Some multivariate Chebyshev inequalities with extensions to continuous parameter processes. *Ann. Math. Statist.* **32** 687–703. MR0148106

[11] CHEN, H., GUO, L. and GAO, A. (1988). Convergence and robustness of the Robbins–Monro algorithm truncated at randomly varying bounds. *Stochastic Process. Appl.* **27** 217–231. MR0931029

[12] CHEN, H. and ZHU, Y.-M. (1986). Stochastic approximation procedures with randomly varying truncations. *Sci. Sinica Ser. A* **29** 914–926. MR0869196

[13] DUFLO, M. (1997). *Random Iterative Systems.* Springer, New York. MR1485774

[14] DVORETZKY, A. (1972). Asymptotic normality for sums of dependent random variables. *Proceedings of the Sixth Berkeley Symposium on Mathematical Statistics and Probability II* 513–535. Univ. California Press, Berkeley. MR0415728

[15] GASEMYR, J. (2003). On an adaptive version of the Metropolis–Hastings with independent proposal distribution. *Scand. J. Statist.* **30** 159–173. MR1965100

[16] GELMAN, A., ROBERTS, G. and GILKS, W. (1995). Efficient Metropolis jumping rules. In *Bayesian Statistics 5* (J. O. Berger, J. M. Bernardo, A. P. Dawid and A. F. M. Smith, eds.) 599–608. Oxford University Press, New York.

[17] GELMAN, A. and RUBIN, D. B. (1992). Inference from iterative simulation using multiple sequences. *Statist. Sci.* **7** 473–483.

[18] GLYNN, P. W. and MEYN, S. P. (1996). A Liapounov bound for solutions of the Poisson equation. *Ann. Probab.* **24** 916–931. MR1404536

[19] HAARIO, H., SAKSMAN, E. and TAMMINEN, J. (2001). An adaptive Metropolis algorithm. *Bernoulli* **7** 223–242. MR1828504

[20] HALL, P. and HEYDE, C. (1980). *Martingale Limit Theory and Its Application*. Academic Press, New York. MR0624435

[21] HOLDEN, L. (1998). Adaptive chains. Technical report, Norwegian Computing Center.

[22] JARNER, S. and HANSEN, E. (2000). Geometric ergodicity of Metropolis algorithms. *Stochastic Process. Appl.* **85** 341–361. MR1731030

[23] KUSHNER, H. and YIN, G. (1997). *Stochastic Approximation Algorithms and Applications*. Springer, New York. MR1453116

[24] MCLEISH, D. (1975). A maximal inequality and dependent strong laws. *Ann. Probab.* **3** 829–839. MR0400382

[25] MENG, X. L. and VAN DYK, D. (1997). The EM algorithm—an old folk song sung to a fast new tune. *J. R. Statist. Soc. Ser B Stat. Methodol.* **59** 511–567. MR1452025

[26] MENGERSEN, K. and TWEEDIE, R. (1996). Rates of convergence of the Hastings and Metropolis algorithms. *Ann. Statist.* **24** 101–121. MR1389882

[27] METROPOLIS, N., ROSENBLUTH, A., ROSENBLUTH, M., TELLER, A. and TELLER, M. (1953). Equations of state calculations by fast computing machines. *J. Chemical Physics* **21** 1087–1091.

[28] MEYN, S. and TWEEDIE, R. (1993). *Markov Chains and Stochastic Stability.* Springer, London. MR1287609

[29] MEYN, S. and TWEEDIE, R. (1994). Computable bounds for convergence rates of Markov chains. *Ann. Appl. Probab.* **4** 981–1011. MR1304770

[30] NUMMELIN, E. (1991). On the Poisson equation in the potential theory of a single kernel. *Math. Scand.* **68** 59–82. MR1124820

[31] ROBERTS, G. and TWEEDIE, R. (1996). Geometric convergence and central limit theorem for multidimensional Hastings and Metropolis algorithms. *Biometrika* **83** 95–110. MR1399158

[32] TITTERINGTON, D., SMITH, A. F. M. and MAKOV, U. (1985). *Statistical Analysis of Finite Mixture Distributions*. Wiley, Chichester. MR0838090

[33] WU, C. (1983). On the convergence properties of the EM algorithm. *Ann. Statist.* **11** 95–103. MR0684867

SCHOOL OF MATHEMATICS  
UNIVERSITY OF BRISTOL  
UNIVERSITY WALK  
BS8 1TW  
UNITED KINGDOM  
E-MAIL: c.andrieu@bristol.ac.uk  
URL: http://www.maths.bris.ac.uk/~maxca/

ENST  
UMR 5061  
46 RUE BARRAULT  
75634 PARIS CEDEX 13  
FRANCE  
E-MAIL: moulines@tsi.enst.fr  
URL: http://www.tsi.enst.fr/~moulines/